%% file: arxiv.tex
\documentclass[a4paper,journal,twoside,web]{ieeecolor}
\usepackage{generic}

\input{doc/header}

\renewcommand{\IEEEPARstart}[2]{#1#2}

\usepackage{fancyhdr}
\fancypagestyle{firstpage}{%
    \fancyfoot[C]{\small
        © 2023 IEEE.
        Personal use of this material is permitted.
        Permission from IEEE must be obtained for all other uses, in any current or future media, including reprinting/republishing this material for advertising or promotional purposes, creating new collective works, for resale or redistribution to servers or lists, or reuse of any copyrighted component of this work in other works.%
    }%
}

\title{A Decomposition Approach to Multi-Agent Systems with Bernoulli Packet Loss}
\author{Christian~Hespe, Hamideh~Saadabadi, Adwait~Datar, Herbert~Werner, and Yang~Tang
\thanks{Christian Hespe, Hamideh Saadabadi, Adwait Datar and Herbert Werner are with the Hamburg University of Technology, Institute of Control Systems, 21073 Hamburg, Germany. (e-mail: \{christian.hespe, hamideh.saadabadi, adwait.datar, h.werner\}@tuhh.de)}%
\thanks{Yang Tang is with the East China University of Science and Technology, Shanghai 200237, China. (e-mail: yangtang@ecust.edu.cn)}}

\begin{document}

\maketitle
\thispagestyle{firstpage}

\begin{abstract}
    \input{doc/abstract}
\end{abstract}

\pagestyle{fancy}

\input{doc/content}
\appendices
\input{doc/appendix}

\bibliographystyle{IEEETran}
\bibliography{arxiv}
\end{document}

%% file: doc/header.tex
\usepackage[T1]{fontenc}
\usepackage[english]{babel}

\usepackage{microtype}

\usepackage{amsmath, amssymb}
\usepackage{amsthm}
\usepackage{mathtools}
\usepackage{commath}
\usepackage[casesstyle, subnum]{cases}
\usepackage{siunitx}

\usepackage{accents}

\usepackage{csquotes}
\usepackage{hyperref}
\usepackage{cite}

\usepackage{booktabs}
\usepackage{enumerate}

\usepackage{tikz}
\usepackage{pgfplots}
\usepackage[caption=false,font=footnotesize]{subfig}

\usetikzlibrary{
 	arrows.meta, 
 	calc,        
 	fit,         
  	positioning, 
}

\pgfplotsset{
    compat=1.16,
    cycle list name=linestyles*,
    every axis plot/.append style={
        thick
    },
    log x ticks with fixed point/.style={
        xticklabel={
            \pgfkeys{/pgf/fpu=true}
            \pgfmathparse{exp(\tick)}%
            \pgfmathprintnumber[fixed relative, precision=3]{\pgfmathresult}
            \pgfkeys{/pgf/fpu=false}
        }
    },
    log y ticks with fixed point/.style={
        yticklabel={
            \pgfkeys{/pgf/fpu=true}
            \pgfmathparse{exp(\tick)}%
            \pgfmathprintnumber[fixed relative, precision=3]{\pgfmathresult}
            \pgfkeys{/pgf/fpu=false}
        }
    }
}

\tikzstyle{agent}        = [draw, circle, fill=gray!30]
\tikzstyle{hidden}       = [dashed, draw=black!80]
\tikzstyle{hidden agent} = [agent, hidden, fill=gray!15]
\tikzstyle{hidden link}  = [link, hidden]
\tikzstyle{link}         = [>=Latex]

\newcommand{\ubar}[1]{\underaccent{\bar}{#1}}

\newcommand{\bbma}{\begin{bmatrix}}
\newcommand{\ebma}{\end{bmatrix}}
\newcommand{\real}{\mathbb{R}}

\DeclareMathOperator{\expect}{\mathbb{E}}
\DeclareMathOperator{\prob}{Pr}
\DeclareMathOperator{\trace}{tr}

\newtheorem{theorem}{Theorem}
\newtheorem{lemma}[theorem]{Lemma}
\newtheorem{corollary}[theorem]{Corollary}

\theoremstyle{definition}
\newtheorem{assumption}{Assumption}
\newtheorem{definition}{Definition}

\theoremstyle{remark}
\newtheorem*{remark}{Remark}

%% file: doc/abstract.tex
In this paper, we extend the decomposable systems framework to multi-agent systems with Bernoulli distributed packet loss with uniform probability.
The proposed sufficient analysis conditions for mean-square stability and \textit{H\textsubscript{2}}-performance -- which are expressed in the form of linear matrix inequalities -- scale linearly with increased network size and thus allow to analyse even very large-scale multi-agent systems.
A numerical example demonstrates the potential of the approach by application to a first-order consensus problem.

%% file: doc/content.tex
\section{Introduction}\label{sec:intro}
\IEEEPARstart{C}{ontrolling} large-scale networks of dynamic systems is a challenging problem that has attracted a lot of research interest.
Due to their vast size, systematic centralized controller synthesis or system analysis quickly become infeasible due to computational demands.
For that reason, decentralized or distributed approaches have become the methods of choice for this class of systems \cite{Massioni2009}.

One particular type of such large-scale systems are multi-agent systems (MAS), in which groups of simple systems -- called agents -- collectively solve tasks by applying agent-level rules.
Examples for such tasks include formation control, distributed estimation or source seeking \cite{Mesbahi2010}.
In order to analyse these MAS, the flexible and powerful framework of \emph{decomposable systems} introduced by Massioni and Verhaegen in \cite{Massioni2009} can be employed.
It is built upon the idea of decoupling the MAS into smaller \emph{modal subsystems} and analysing these subsystems independently, a technique which was originally introduced in \cite{Fax2004} for stability analysis only.
By decoupling the analysis, the framework improves the scalability in terms of computational complexity from quadratic to linear, in some instances even constant, in the number of agents \cite{Massioni2009}.
Originally proposed for linear-time invariant (LTI) systems, the framework has been extended to linear parameter-varying systems \cite{Hoffmann2013} and analysis using integral quadratic constraints \cite{Eichler2013}.

An important aspect of MAS is how the exchange of information is implemented.
Depending on the requirements, relative measurements or a communication network are preferable.
In this paper, we will be focusing on the latter and consider the case where the communication network is subject to stochastic uncertainty in form of lost information.
More specifically, we investigate how to analyse the effect of packet loss described by independent Bernoulli distributed random variables with uniform probability on stability and performance in a scalable manner.
As noted by the authors of \cite{Ma2020}, most existing work on networked MAS with stochastic packet loss assumes identical loss, i.e., that all communication links fail at the same time, an assumption very few systems satisfy in practice.
Amongst others, this scenario is studied in \cite{Zhang2017, Wang2018} and \cite{Xu2020} for Bernoulli and Markov packet loss models, respectively.
On the other hand, there are approaches that consider not identical \emph{loss} but uniform \emph{loss probability}, e.g. \cite{Mesbahi2010, Patterson2010, Wu2012, Ghadami2012}.
All four assume \emph{symmetric} loss, i.e., that link failure is identical in both directions.
Finally, Bernoulli packet loss with non-uniform probabilities and independent links is considered in \cite{Zhang2012} for directed tree graphs using only a lower bound on the transmission probabilities and in \cite{Ma2020} for general graphs with know probability for each link.

Of the aforementioned papers, only \cite{Mesbahi2010, Patterson2010, Ghadami2012} consider system performance in addition to stability, the first two in terms of the convergence rate, the third using the \(L_2\) system norm.
Another important performance measure for MAS is the \(H_2\)-norm, see \cite{Massioni2009, Ghadami2012, Stoorvogel2019} and \cite{Raza2022} amongst others.
A stochastic generalization of this norm for Markov jump linear systems (MJLS) was introduced for optimal control in \cite{Costa1997} and used for optimal filtering in \cite{Fioravanti2008}.
An existing approach for analysing large MAS with MJLS can be found in \cite{Lee2015}.
However, while the conditions scale linearly with the number of agents, they scale exponentially with the maximum vertex degree and are thus intractable for many MAS.

Modelling packet loss with identically Bernoulli distributed random variables is invalid in many real-world scenarios.
Nonetheless, this paper provides a first step towards \emph{scalable} analysis of MAS with more realistic networking models.

\subsection{Contributions}\label{sec:intro_contribution}
The main contribution of this paper are the sufficient analysis conditions for mean-square stability and \(H_2\)-performance of MAS presented in Theorems~\ref{thm:ms_stability_decoupled} and \ref{thm:mjls_h2_decoupled} that scale \emph{linearly} with the number of agents in presence of \emph{non-identical} Bernoulli distributed packet loss with uniform probability.
The conditions are based on extending the decomposable systems framework to stochastic jump linear systems and the analytic calculation of the expected Laplacian matrices in Lemma~\ref{lem:expected_laplacians}.
Similar analytic calculations have been presented before in \cite{Wu2012}, however, in contrast to previous works and at the cost of losing necessity, the current paper does \emph{not} rely on having symmetric packet loss and brings out the inherent structure of the expected Laplacian matrices allowing for decomposition, which is exploited in Lemma~\ref{lem:simultaneous_diagonalization}.
Two further smaller contributions are necessary conditions in Theorem~\ref{thm:necessary_conditions} supporting the sufficient conditions and an analysis approach for uncertain transmission probabilities and communication topologies that is based on convexity arguments.

\subsection{Outline}\label{sec:intro_outline}
Following this introduction, Section~\ref{sec:problem} proceeds with defining notation, setting up the problem and extending the decomposable systems framework.
Section~\ref{sec:laplacian} contains the calculation of the expected Laplacians.
The main results are presented in Section~\ref{sec:scalable}, followed by a numerical example in Section~\ref{sec:example}.
Finally, the paper is concluded in Section~\ref{sec:conclusions}.

\section{Problem Statement}\label{sec:problem}
\subsection{Notation and Definitions}\label{sec:problem_notation}
We let \(I_N\) denote the \(N \times N\) identity matrix and \(\mathbf{1}_N\) the vector in \(\real^N\) with all entries equal to 1.
\(M \succ (\succeq)\ 0\) or \(M \prec (\preceq)\ 0\) mean that \(M\) is positive or negative (semi-) definite.
\(M_1 \otimes M_2\) is the Kronecker product, which has the mixed product property \((A \otimes B)(C \otimes D) = (AC) \otimes (BD)\).
\(\mathbb{I}_S\) denotes the set-membership indicator function defined as
\begin{equation}
    \mathbb{I}_S(x) \coloneqq \begin{cases}
        1 & if \(x \in S\), \\
        0 & else.
    \end{cases}
\end{equation}
Depending on the context, we use \(\|z\|\) either for the Euclidean vector norm, the induced matrix 2-norm, or the 2-norm for (stochastic) signals defined by \(\|z\|^2 \coloneqq \sum_{k = 0}^\infty \expect\left[z^T(k) z(k)\right]\).

The interconnections between agents are modelled using graphs~\(\mathcal{G} \coloneqq (\mathcal{V}, \mathcal{E})\), which are composed of the vertex set \(\mathcal{V} = \{v_1, v_2, \ldots, v_N\}\) and the edge set \(\mathcal{E} \subset \mathcal{V} \times \mathcal{V}\), where an edge \(e_{ij} \coloneqq (v_j, v_i)\) is read as pointing from \(v_j\) to \(v_i\) and \(e_{ii} \notin \mathcal{E}\).
\(\mathcal{G}\) is called undirected if \(e_{ij} \in \mathcal{E} \Leftrightarrow e_{ji} \in \mathcal{E}\).
The set \(\mathcal{N}_i^- \coloneqq \{v_j \in \mathcal{V} : e_{ij} \in \mathcal{E}\}\) is called the in-neighbourhood of \(v_i\) and its cardinality \(d_i^- \coloneqq |\mathcal{N}_i^-|\) is the in-degree of \(v_i\).
Equivalently, define the out-neighbourhood \(\mathcal{N}_i^+ \coloneqq \{v_j \in \mathcal{V} : e_{ji} \in \mathcal{E}\}\) and out-degree \(d_i^+ \coloneqq |\mathcal{N}_i^+|\).
If for every vertex in \(\mathcal{V}\) the in- and out-degree are identical, \(\mathcal{G}\) is said to be balanced.
A sequence of vertices is called directed path on \(\mathcal{G}\) if \(e_{ji} \in \mathcal{E}\) for all pairs of consecutive vertices \((v_i, v_j)\).
If there exists a directed path from all \(v_r \in \mathcal{V}\) to all other \(v_i \in \mathcal{V} \setminus {v_r}\), \(\mathcal{G}\) is said to be strongly-connected.
The transpose \(\mathcal{G}^T\) is defined as the graph in which the direction of every edge is inverted, i.e. \(\mathcal{G}^T \coloneqq (\mathcal{V}, \mathcal{E}^T)\) with \(e_{ij} \in \mathcal{E} \Leftrightarrow e_{ji} \in \mathcal{E}^T\).

For a graph~\(\mathcal{G}\), define element-wise the Laplacian matrix \(L(\mathcal{G}) \coloneqq [l_{ij}(\mathcal{G})]\), where
\begin{equation}\label{eq:laplacian_nominal}
    l_{ij}(\mathcal{G}) \coloneqq \begin{cases}
        -1     & if \(i \neq j\) and \(v_j \in    \mathcal{N}_i^-\), \\
         0     & if \(i \neq j\) and \(v_j \notin \mathcal{N}_i^-\), \\
         d_i^- & if \(i = j\).
    \end{cases}
\end{equation}
\(L(\mathcal{G})\) is symmetric if and only if \(\mathcal{G}\) is undirected.
We will drop the argument from the notation if the corresponding graph can be determined from context.

\subsection{Jump Linear Systems for Modelling Packet Loss}\label{sec:problem_mjls}
The focus of this paper are MAS which are subject to stochastic packet loss.
This kind of system cannot be modelled in a time-invariant manner, since loss of packets means that connections between individual agents break momentarily and thus the interconnection topology between agents is time-varying.
For this reason, we will use a special case of MJLS to model the MAS.

An MJLS is a discrete-time, switched linear system whose switching is controlled by a corresponding Markov chain.
At every time instance, the MJLS is in exactly one of \(m\) possible modes, where each mode can have a different dynamic behaviour.
It is described by the state-space system
\begin{equation}\label{eq:mjls}
    G: \,\left\{\enspace\begin{aligned}
        x(k+1) &= A_{\sigma(k)} x(k) + B_{\sigma(k)} w(k) \\
        z(k)   &= C_{\sigma(k)} x(k) + D_{\sigma(k)} w(k),
    \end{aligned}\right.
\end{equation}
where \(x(k) \in \real^{N n_x}\) is the dynamic state, \(\sigma(k) \in \mathcal{K} \coloneqq \left\{i \in \mathbb{N} : 1 \leq i \leq m\right\}\) is the state of the Markov chain and \(w(k) \in \real^{N n_w}\) and \(z(k) \in \real^{N n_z}\) are the performance input and output, respectively.
\(N\) denotes the number of agents in the system, the initial state of the system is \(x(0) = x_0\) and the Markov chain is initially distributed according to \(\sigma(0) = \sigma_0\).
For each mode \(i \in \mathcal{K}\), the dynamics of the system are governed by the matrices \(A_i\), \(B_i\), \(C_i\) and \(D_i\).
Note that system~\eqref{eq:mjls} does not have a control input or measured output since this paper is only concerned with system analysis in contrast to controller synthesis.
\eqref{eq:mjls} should thus be considered as a closed-loop model, containing an agent model and potentially a controller.

In this paper, we only consider the case where the switching probability of the Markov chain is independent of the chain's state, thus the distribution of \(\{\sigma(k)\}\) is stationary and described by
\begin{equation}\label{eq:markov_chain}
    \prob\big(\sigma(k) = i\big) = t_i,
\end{equation}
for all \(k \geq 0\).

There is a variety of definitions for stability in the context of MJLS.
Amongst them are stability in expectation, almost sure stability and mean-square stability (MSS).
Here we will focus on the latter.
In comparison, MSS has the advantage that it is easy to test for and implies stability as in the other two definitions \cite{Costa2005}.

\begin{definition}[Mean-Square Stability \cite{Costa2005}]\label{def:ms_stability}
    The MJLS \eqref{eq:mjls} is mean-square stable if
    \begin{align*}
        \lim\limits_{k \to \infty} \expect\left[\|x(k)\|\right] = 0 & &
        \text{and} & &
        \lim\limits_{k \to \infty} \expect\left[\|x(k) x^T(k)\|\right] = 0
    \end{align*}
    for all initial conditions~\(x_0\) and initial distributions~\(\sigma_0\).
\end{definition}

In the following, we will often refer to an MJLS as stable if it is MSS.
As shown in \cite{Costa2005}, stability of the individual modes of an MJLS is neither necessary nor sufficient for MSS.
Instead, we will make use of the following linear matrix inequality (LMI) based stability test:

\begin{theorem}[LMI Condition for MSS \cite{Costa1993}]\label{thm:ms_stability}
    The MJLS \eqref{eq:mjls} is mean-square stable if and only if there exists a \(Q \succ 0\) such that
    \begin{equation}\label{eq:ms_stability_lmi}
        \sum_{i \in \mathcal{K}} t_{i} A_i^T Q A_i - Q \prec 0.
    \end{equation}
\end{theorem}
\begin{remark}
    Note that we can express the above LMI using an unconditional expectation.
    Thus, \eqref{eq:ms_stability_lmi} is equivalent to
    \begin{equation}\label{eq:ms_stability_lmi_expectation}
        \expect\left[A_\sigma^T Q A_\sigma\right] - Q \prec 0,
    \end{equation}
    where the expectation is taken with respect to \(\sigma\).
\end{remark}

This theorem is a specialization of the general stability test from \cite{Costa1993} to MJLS with state-independent switching probabilities as given in \eqref{eq:markov_chain}.
Compared to the general case, this theorem results in a sizeable reduction in computational complexity, since the stability test contains only a single matrix variable \(Q\) and a single LMI constraint, instead of having one of both for each mode.
Still, it is necessary to enumerate all modes in \eqref{eq:ms_stability_lmi}, which renders the analysis of systems with numerous agents intractable.
For the specific system structure that is introduced in the next subsection, we will develop an approach that eliminates the need for mode enumeration.

In addition to MSS, we consider system performance in terms of the \(H_2\)-norm from input~\(w\) to output~\(z\).
For the special case of mode independent transition probabilities in the Markov chain, the norm is defined as follows:

\begin{definition}[\(H_2\)-norm for MJLS \cite{Fioravanti2013}]\label{def:mjls_h2}
    The \(H_2\)-norm of the stable MJLS~\eqref{eq:mjls} is defined as
    \begin{equation*}
        \|G\|_{H_2}^2 \coloneqq \sum_{i \in \mathcal{K}} \sum_{s = 1}^{n_w} t_i \|z^{s,i}\|^2,
    \end{equation*}
    where \(z^{s,i}\) is the response of \(G\) to a discrete impulse applied into the \(s\)th input with \(x_0 = 0\) and \(\sigma_0 = i\).
\end{definition}

Similar to Theorem~\ref{thm:ms_stability} for MSS, we can exploit the stationarity of the transition probabilities to obtain an analysis condition in two variables and LMI constraints.
For general MJLS, the corresponding condition requires two LMIs and variables \emph{for each mode}, resulting in much larger computational cost.
The procedure to obtain this simplified analysis condition was introduced in \cite{Fioravanti2012}.

\begin{theorem}[LMI condition for MJLS \(H_2\)-norm \cite{Fioravanti2013}]\label{thm:mjls_h2}
    Given the stable MJLS~\eqref{eq:mjls}, \(\|G\|_{H_2} < \gamma\) if and only if there exist a \(Q \succ 0\) and a symmetric \(Z\) with \(\trace\left(Z\right) < \gamma^2\) such that%
    \begin{subequations}\label{eq:mjls_h2_lmi}\begin{align}%
        \sum_{i \in \mathcal{K}} t_i \left(A_i^T Q A_i + C_i^T C_i\right) - Q &\prec 0, \label{eq:mjls_h2_gramian_lmi} \\
        \sum_{i \in \mathcal{K}} t_i \left(B_i^T Q B_i + D_i^T D_i\right) - Z &\prec 0. \label{eq:mjls_h2_trace_lmi}
    \end{align}\end{subequations}
\end{theorem}
\begin{remark}
    Theorem~\ref{thm:mjls_h2} reformulates the result from \cite{Fioravanti2013} by introducing \(Z\).
    To see that both are equivalent, use \(Y - X \prec 0 \Rightarrow \trace\left(X\right) > \trace\left(Y\right)\) with symmetric \(X,Y\) for the first direction and chose \(Z = \varepsilon I + \sum_{i \in \mathcal{K}} t_i (B_i^T Q B_i + D_i^T D_i)\) with sufficiently small \(\varepsilon > 0\) for the other.
\end{remark}

For the same reason as for Theorem~\ref{thm:ms_stability}, direct application of the above result to large MAS would quickly lead to numerically intractable problems.
We introduce a subset of jump systems for which the computational complexity can be vastly reduced next.

\subsection{Decomposable Jump Linear Systems}\label{sec:problem_decomposable}
Coming from the general MJLS in subsection~\ref{sec:problem_mjls}, this paper considers systems with a specific structure in their state-space matrices which allows to utilize the decomposable systems framework introduced by Massioni and Verhaegen in \cite{Massioni2009}.
According to their definition, a matrix \(M\) is said to be decomposable if it can be split up into a decoupled component \(M^d\) and a coupled component \(M^c\) as \(M = I_N \otimes M^d + P \otimes M^c\), where \(P\) is called the pattern matrix.
Moreover, an LTI system is called decomposable if all matrices of its state-space representation are decomposable with respect to the \emph{same} pattern matrix.

Applying this concept to the MJLS~\eqref{eq:mjls} means that \(A_i\), \(B_i\), \(C_i\) and \(D_i\) must be decomposable, we do however \emph{not} insist on having the same pattern matrix for all modes \(i \in \mathcal{K}\).
On the contrary, we will assume the pattern matrix is the only part of the system that changes between the modes.
This choice is motivated by the fact that -- in the context of networked multi-agent systems -- the pattern matrix is given by the graph Laplacian and that the communication graph is a stochastic process due to packet loss.
We then introduce the nominal graph \(\mathcal{G}^0 = (\mathcal{V}, \mathcal{E}^0)\) and its corresponding Laplacian \(L^0 := L(\mathcal{G}^0)\).
All together, this leads to the decomposable MJLS
\begin{equation}\label{eq:mjls_decomposable}
    \hat{G}: \,\left\{\enspace\begin{aligned}
        x(k+1) &= \left(I_N \otimes A^d + L\big(\mathcal{G}_{\sigma(k)}\big) \otimes A^c\right) x(k) \\
               +& \left(I_N \otimes B^d + L\big(\mathcal{G}_{\sigma(k)}\big) \otimes B^c\right) w(k), \\
        z(k)   &= \left(I_N \otimes C^d + L\big(\mathcal{G}_{\sigma(k)}\big) \otimes C^c\right) x(k) \\
               +& \left(I_N \otimes D^d + L\big(\mathcal{G}_{\sigma(k)}\big) \otimes D^c\right) w(k),
    \end{aligned}\right.
\end{equation}
where \(\mathcal{G}_i \coloneqq (\mathcal{V}, \mathcal{E}_i)\) and where \(\mathcal{E}_i \subseteq \mathcal{E}^0\) is the subset of edges that successfully transmit a packet in mode \(i\) of the MJLS.
Analogously to \(L^0\), define \(L_i \coloneqq L\left(\mathcal{G}_i\right)\) as shorthand notation.

More specifically, consider a stationary stochastic process \(\{\alpha_{ij}(k)\}\) for each \(e_{ij} \in \mathcal{E}^0\), where \(\alpha_{ij}(k) \in \{0, 1\}\).
Here, \(\alpha_{ij}(k) = 1\) means the edge \(e_{ij}\) is active, or equivalently that the packet is transmitted, while \(\alpha_{ij}(k) = 0\) means \(e_{ij}\) is inactive and the packet is lost.
The edges might fail asymmetrically, i.e., it might happen that \(\alpha_{ij}(k) \neq \alpha_{ji}(k)\).
In the following, we assume that the stochastic processes are Bernoulli distributed and independent in time.
Furthermore, at any given time instant, the packet loss between two different pairs of vertices is assumed to be independent.
This is formalized in the following assumption.

\begin{assumption}\label{ass:bernoulli}
    The stochastic processes \(\{\alpha_{ij}(k)\}\) are \emph{partially} independent and identically Bernoulli distributed such that, for all \(k, k' \geq 0\), \(e_{ij}, e_{rs} \in \mathcal{E}^0\), we have
    \begin{align}\label{eq:bernoulli_probability}
        \prob\left(\alpha_{ij}(k) = 1\right) &= p & \prob\left(\alpha_{ij}(k) = 0\right) &= 1 - p
    \end{align}
    with \(p \in [0, 1]\) and \(\alpha_{ij}(k)\) and \(\alpha_{rs}(k')\) are independent random variables whenever \(k \neq k'\) or \((r,s) \neq (i,j) \wedge (r,s) \neq (j,i)\).
\end{assumption}
\begin{remark}
    For many real world scenarios, modelling packet loss as independent Bernoulli distributed random variables with uniform probability is an idealization. Similar to \cite{Schenato2007}, we proceed in this way for reasons of mathematical tractability.
    Note that compared to assuming identical or symmetric loss as e.g.\@ in \cite{Zhang2017, Xu2020, Wu2012}, Assumption~\ref{ass:bernoulli} is closer to reality due to allowing opposing links to be correlated or not.
\end{remark}

To map from the stochastic processes \(\{\alpha_{ij}(k)\}\) to the MJLS~\eqref{eq:mjls_decomposable}, define a function \(\nu : \mathcal{E}^0 \to \left\{1, \ldots, |\mathcal{E}^0|\right\}\) that assigns each \(e_{ij} \in \mathcal{E}^0\) a unique integer.
Then, we have
\begin{equation}\label{eq:sigma_mapping}
    \sigma(k) = 1 + \sum\limits_{e_{ij} \in \mathcal{E}^0} \alpha_{ij}(k) 2^{\nu\left(e_{ij}\right) - 1}
\end{equation}
and accordingly \(m = 2^{|\mathcal{E}^0|}\).
The map from \(\alpha_{ij}(k)\) to \(\sigma(k)\) is bijective, such that we can equivalently represent the Bernoulli packet loss model in form of the MJLS.
Thus, our communication model has two parameters:
The graph \(\mathcal{G}^0\) and the probability of successful transmission \(p\).

To reap maximum benefit from introducing the decomposable system framework, we will impose that the matrix variable \(Q\) has block repeated structure.
While this may be a conservative choice, it allows generating stability and performance tests that are particularly easy to check.

\begin{corollary}[MSS for Decomposable Jump Systems]\label{cor:ms_stability_decomposable}
    The decomposable jump system \eqref{eq:mjls_decomposable} is mean square stable if there exists a \(Q \succ 0\) such that
    \begin{equation}\label{eq:ms_stability_decomposable_lmi}\begin{aligned}
        I_N \otimes \big(A^{dT} Q A^d - Q\big) + \expect\big[L_\sigma\big] \otimes \big(A^{dT} Q A^c\big) +&\\
        \expect\big[L_\sigma^T\big] \otimes \big(A^{cT} Q A^d\big) + \expect\big[L_\sigma^T L_\sigma\big] \otimes \big(A^{cT} Q A^c\big) &\prec 0.
    \end{aligned}\end{equation}
\end{corollary}
\begin{proof}
    Take the LMI condition \eqref{eq:ms_stability_lmi_expectation} and insert \(Q = I_N \otimes \tilde{Q}\) as well as the decomposable MJLS from \eqref{eq:mjls_decomposable} for \(A_\sigma\), resulting in
    \begin{align*}
        \expect\Big[(I_N \otimes A^d + L_\sigma \otimes A^c)^T (I_N \otimes \tilde{Q})&\\
        (I_N \otimes A^d + L_\sigma \otimes &A^c)\Big] - I_N \otimes \tilde{Q} \prec 0.
    \end{align*}
    Using the mixed product rule and the commutation property \((M_1 \otimes I)(I \otimes M_2) = (I \otimes M_2)(M_1 \otimes I)\) we obtain
    \begin{align*}
        I_N \otimes \big(A^{dT} \tilde{Q} A^d - \tilde{Q}\big) + \expect\Big[L_\sigma \otimes \big(A^{dT} \tilde{Q} A^c\big) +& \\
        L_\sigma^T \otimes \big(A^{cT} \tilde{Q} A^d\big) + \big(L_\sigma^T L_\sigma\big) \otimes \big(A^{cT} \tilde{Q} &A^c\big)\Big] \prec 0,
    \end{align*}
    from where we get to \eqref{eq:ms_stability_decomposable_lmi} by linearity of the expectation, \(\expect[X \otimes M] = \expect[X] \otimes M\), which holds if \(X\) is a random matrix and \(M\) a constant, and renaming \(\tilde{Q} \to Q\).
\end{proof}
\begin{remark}
    The only source of conservatism in Corollary~\ref{cor:ms_stability_decomposable} is the assumption on \(Q\) to have block repeated structure.
    A similar result can be obtained without this assumption, however, isolating the expectation of \(L_\sigma^T L_\sigma\) would not be possible, since the commutation property cannot be used.
    Instead, one would have to consider a weighted squared expectation of the form \(\expect[(L_\sigma^T \otimes I) Q (L_\sigma \otimes I)]\), similar to \cite[Lemma 1]{Wu2012} but with additional Kronecker products.
    In that case, the analysis conditions cannot be decomposed using the approach proposed in Section~\ref{sec:scalable_decomposition} below.
\end{remark}

A similar corollary can be derived for the \(H_2\)-performance conditions in Theorem~\ref{thm:mjls_h2}.
However, since the derivation would be analogous to the proof of Corollary~\ref{cor:ms_stability_decomposable}, we skip this intermediate result.

\section{Expected Laplacian Matrices}\label{sec:laplacian}
From Corollary~\ref{cor:ms_stability_decomposable}, we have seen how the expectation of the Laplacian is essential in determining if a decomposable MJLS is MSS or not.
We thus derive an analytic calculation of the expectation in terms of \(\mathcal{G}^0\) and \(p\) in the following.

As preparation, notice how the elements of the Laplacian change compared to \eqref{eq:laplacian_nominal} when packet loss is introduced.
Using the element-wise notation \(L_{\sigma(k)} = [l_{ij}^\sigma(k)]\), the stochastic Laplacian is given by
\begin{equation}\label{eq:random_laplacian}
    l_{ij}^\sigma(k) = \begin{cases}
        -\alpha_{ij}(k)                                  & if \(i \neq j\) and \(v_j \in \mathcal{N}_i^-\), \\
        0                                                & if \(i \neq j\) and \(v_j \notin \mathcal{N}_i^-\), \\
        \sum\limits_{v_m \in \mathcal{N}_i^-} \alpha_{im}(k) & if \(i = j\).
    \end{cases}
\end{equation}
As noted above, \(\{L_{\sigma(k)}\}\) is a stochastic process due to packet loss.
Since \(\{\alpha_{ij}(k)\}\) and \(\{L_{\sigma(k)}\}\) are stationary by Assumption~\ref{ass:bernoulli}, we will drop the index \(k\) in the remainder of the paper when referring to an instance of these processes.

\begin{lemma}[Expected Laplacian Matrices]\label{lem:expected_laplacians}
    Given the nominal graph~\(\mathcal{G}^0\) and packet loss according to Assumption~\ref{ass:bernoulli}, we have
    \begin{align*}
        \expect[L_\sigma]            &= p L(\mathcal{G}^0), \\
        \expect[L_\sigma^T L_\sigma] &= p^2 L(\mathcal{G}^0)^T L(\mathcal{G}^0) + p(1-p) \big(L(\mathcal{G}^0) + L(\mathcal{G}^{0T})\big).
    \end{align*}
\end{lemma}
\begin{proof}
    See Appendix~\ref{app:lemma_expectation}.
\end{proof}

Lemma~\ref{lem:expected_laplacians} enables us to calculate the expected Laplacians analytically from the two parameters \(\mathcal{G}^0\) and \(p\), which allows applying Corollary~\ref{cor:ms_stability_decomposable} and Theorem~\ref{thm:mjls_h2} effectively without expensive numerical calculation of the expectations by enumeration of all modes.
Note that \(L(\mathcal{G}^T) = L(\mathcal{G})^T\) if and only if \(\mathcal{G}\) is balanced.
For certain graphs \(\mathcal{G}^0\), we can further exploit the following diagonalizability property:

\begin{lemma}[Simultaneous Diagonalizability]\label{lem:simultaneous_diagonalization}
    Given the nominal graph \(\mathcal{G}^0\) and packet loss according to Assumption~\ref{ass:bernoulli}, there exists a similarity transformation~\(U\) that diagonalizes \(\expect[L_\sigma]\), \(\expect[L_\sigma^T]\) and \(\expect[L_\sigma^T L_\sigma]\) if and only if \(L(\mathcal{G}^0)\) is normal, i.e. \(L(\mathcal{G}^0) L(\mathcal{G}^0)^T = L(\mathcal{G}^0)^T L(\mathcal{G}^0)\).
\end{lemma}
\begin{proof}
    According to \cite[p. 62]{Horn2012}, there exists a similarity transformation that diagonalizes two diagonalizable matrices at the same time if and only if the matrices commute.
    As shown in \cite{Wu1995}, \(L(\mathcal{G}^0)\) being normal and having zero row sum implies that it has zero column sum as well, therefore \(\mathcal{G}^0\) is balanced and \(L(\mathcal{G}^{0T}) = L(\mathcal{G}^0)^T\).
    By Lemma~\ref{lem:expected_laplacians}, we thus need to show that \(L^0\), \(L^{0T}\) and \(L^{0T} L^0\) commute.
    From the definition of normality, this is trivial for the first pair and easy to verify for \(L^0\) and \(L^{0T} L^0\).
    Conversely, if there exists a transformation that diagonalizes both \(\expect[L_\sigma]\) and \(\expect[L_\sigma^T]\), then \(L^0\) and \(L^{0T}\) commute, implying that \(L^0\) is normal.
\end{proof}

For some scenarios in the context of MAS control and distributed consensus, normality of the Laplacian is too restrictive for Lemma~\ref{lem:simultaneous_diagonalization} to be applicable.
In particular, leader-follower schemes cannot be handled, since they require unbalanced communication graphs.

\section{Scalable Analysis with Packet Loss}\label{sec:scalable}
\subsection{Decomposed Analysis LMIs}\label{sec:scalable_decomposition}
With the results from Section~\ref{sec:laplacian}, we can now formulate our final MSS and \(H_2\)-performance analysis conditions for decomposable MJLS.
The motivation for defining a decomposable system like in \cite{Massioni2009} is that we can decouple the system as long as we can diagonalize the pattern matrix.
Assuming there exists a transformation \(U\) such that \(U P U^{-1}\) is diagonal -- with \(P\) being the pattern matrix --, then \((U \otimes I)\) decouples the system matrices.
In particular, if the underlying graph \(\mathcal{G}^0\) is undirected, such a transformation is guaranteed to exist with \(U U^T = I\).
We will thus make the following assumption in the remainder of the paper:

\begin{assumption}\label{ass:undirected_graph}
    The communication graph \(\mathcal{G}^0\) is undirected.
\end{assumption}

Assumption~\ref{ass:undirected_graph} restricts the classes of MAS the following results can be applied to.
Note, however, that it is different from assuming the packet loss is symmetric, which would be equivalent to assuming \emph{all} \(\mathcal{G}_i\) are undirected, in contrast to just \(\mathcal{G}^0\).
Applied to Corollary~\ref{cor:ms_stability_decomposable}, this gives rise to the following stability test consisting of a set of decoupled LMIs:

\begin{theorem}[Decomposed MSS Test]\label{thm:ms_stability_decoupled}
    Given the MJLS~\eqref{eq:mjls_decomposable} with nominal communication graph~\(\mathcal{G}^0\) and packet loss satisfying Assumptions~\ref{ass:bernoulli} and \ref{ass:undirected_graph}, the MJLS is mean-square stable if there exists a \(Q \succ 0\) such that
    \begin{equation}\label{eq:ms_stability_decoupled_lmi}\begin{aligned}
        \big(A^d + p \lambda_i A^c\big)^T Q \big(&A^d + p \lambda_i A^c\big) - Q \\
        &+ 2p(1-p) \lambda_i A^{cT} Q A^c \prec 0
    \end{aligned}\end{equation}
    for all \(i \in \{1, \ldots, N\}\), where \(\lambda_i\) are the eigenvalues of \(L^0\).
\end{theorem}
\begin{proof}
    Since \(\mathcal{G}^0\) is undirected, \(L^0\) is symmetric, and we know from Lemma~\ref{lem:simultaneous_diagonalization} that \(\expect[L_\sigma]\), \(\expect[L_\sigma^T]\) and \(\expect[L_\sigma^T L_\sigma]\) can be diagonalized using an orthogonal matrix \(U\).
    Apply a congruence transformation to \eqref{eq:ms_stability_decomposable_lmi} from Corollary~\ref{cor:ms_stability_decomposable} by multiplying with \((U \otimes I_{n_x})\) and \((U^T \otimes I_{n_x})\) from the left and right, respectively.
    Using the mixed product rule and commutation property of the Kronecker product, this results in
    \begin{align*}
        &I_N \otimes \big(A^{dT} Q A^d - Q\big) + \big(U \expect\big[L_\sigma^T L_\sigma\big] U^T\big) \otimes \big(A^{cT} Q A^c\big) \\
        &\qquad+ \big(U \expect\big[L_\sigma\big] U^T\big) \otimes \big(A^{dT} Q A^c\big) \\
        &\qquad+ \big(U \expect\big[L_\sigma^T\big] U^T\big) \otimes \big(A^{cT} Q A^d\big) \prec 0.
    \end{align*}
    By Lemma~\ref{lem:expected_laplacians}, we have \(U \expect\big[L_\sigma\big] U^T = U \expect\big[L_\sigma^T\big] U^T = p\Lambda\) and \(U \expect\big[L_\sigma^T L_\sigma\big] U^T = p^2\Lambda^2 + 2p(1-p) \Lambda\), where \(\Lambda\) is a diagonal matrix containing the eigenvalues of \(L^0\).
    After the transformation, we have
    \begin{align*}
        I_N \otimes \big(A^{dT} Q A^d - Q\big) + p\Lambda \otimes \big(A^{dT} Q A^c + A^{cT} Q A^d\big)& \\
        + \big(p^2\Lambda^2 + 2p(1-p) \Lambda\big) \otimes \big(A^{cT} Q A^c&\big) \prec 0,
    \end{align*}
    which is a block-diagonal matrix inequality.
    Finally, \eqref{eq:ms_stability_decoupled_lmi} can be obtained by algebraic matrix manipulations and considering the blocks independently.
\end{proof}

Theorem~\ref{thm:ms_stability_decoupled} has multiple advantages in terms of computational complexity compared to the original stability test in Theorem~\ref{thm:ms_stability}.
The first and most impactful is replacing the mode enumeration of the MJLS by the formula given in Lemma~\ref{lem:expected_laplacians}.
Since the number of modes scales at least with \(2^N\) for strongly-connected graphs -- there exists at least one edge per agent --, the original formulation has exponential complexity while the analytic calculation scales quadratically.
The second improvement comes from decomposing the single large constraint on the whole network into multiple smaller ones with the size of a single agent.
In analogy to the modal subsystems from \cite{Massioni2009}, we may term these as \emph{modal constraints}.
Instead of scaling the number of variables and constraints quadratically with the agent count, the decoupled formulation is of constant complexity in the variables and linear complexity in the constraints.
Analogous steps can be applied to the \(H_2\)-performance analysis LMIs from Theorem~\ref{thm:mjls_h2}:

\begin{theorem}[Decomposed \(H_2\)-Performance]\label{thm:mjls_h2_decoupled}
    Given the MJLS~\eqref{eq:mjls_decomposable} with nominal communication graph \(\mathcal{G}^0\) and packet loss satisfying Assumptions~\ref{ass:bernoulli} and \ref{ass:undirected_graph}, \(\hat{G}\) is mean-square stable and \(\|\hat{G}\|_{H_2} < \gamma\) if there exist a \(Q \succ 0\) and symmetric \(Z_i\) with \(\sum_{i = 1}^N \trace\left(Z_i\right) < \gamma^2\) such that%
    \begin{subequations}\label{eq:mjls_h2_decoupled_lmi}\begin{align}\begin{split}%
            \bar{A}_i^T &Q \bar{A}_i + \bar{C}_i^T \bar{C}_i - Q\\
            &+ 2p(1-p) \lambda_i \big(A^{cT} Q A^c + C^{cT} C^c\big) \prec 0 \label{eq:mjls_h2_decoupled_gramian_lmi}
        \end{split}\\\begin{split}
            \bar{B}_i^T &Q \bar{B}_i + \bar{D}_i^T \bar{D}_i - Z_i \\
            &+ 2p(1-p) \lambda_i \big(B^{cT} Q B^c + D^{cT} D^c\big) \prec 0 \label{eq:mjls_h2_decoupled_trace_lmi}
    \end{split}\end{align}\end{subequations}
    for all \(i \in \{1, \ldots, N\}\), where \(\lambda_i\) are the eigenvalues of \(L^0\), \(\bar{A}_i\) denotes \(A^d + p \lambda_i A^c\) and equivalently for \(\bar{B}_i\), \(\bar{C}_i\) and \(\bar{D}_i\).
\end{theorem}
\begin{proof}
    See Appendix~\ref{app:theorem_decomposed_h2}.
\end{proof}

The computational performance improvements achieved by Theorems~\ref{thm:ms_stability_decoupled} and \ref{thm:mjls_h2_decoupled} come at the cost of some conservatism due to imposing that \(Q\) is identical for all modal constraints.
As noted in the remark to Corollary~\ref{cor:ms_stability_decomposable}, this restriction is inherently required to utilize the commutation property of the Kronecker product and thus to apply Lemmas~\ref{lem:expected_laplacians} and \ref{lem:simultaneous_diagonalization} for the calculation of the expected Laplacians.
Calculation of \emph{weighted} expected Laplacians and whether their structure allows for a decomposition of the analysis is subject to further research.
To evaluate how much conservatism is introduced by the restriction to a single \(Q\), we will present a numerical example that demonstrates the trade-off between computational speed and overestimation of the \(H_2\)-norm in Section~\ref{sec:example}.

\subsection{Handling Uncertain Loss Probabilities}\label{sec:scalable_uncertain}
Theorems~\ref{thm:ms_stability_decoupled} and \ref{thm:mjls_h2_decoupled} consider the case where the transmission probability~\(p\) is known exactly.
In practice that is often not the case and only a lower bound \(\ubar{p} \geq 0\) on the transmission probability is known.
If an upper bound \(\bar{p} \leq 1\) is provided in the same vein, whether the MJLS~\eqref{eq:mjls_decomposable} is stable or has \(H_2\)-norm less than \(\gamma\) for a \emph{constant but uncertain} transmission probability can be answered by applying the theorems for all \(p\) in \([\ubar{p}, \bar{p}]\).
However, since \([\ubar{p}, \bar{p}]\) is a real interval, numerical evaluation of the LMI constraints for all such \(p\) is intractable.
Instead, we can make use of the fact that all three LMIs are convex in \(p\) under conditions specified in the following lemma.

\begin{lemma}[Convexity in \(p\)]\label{lem:lmi_convexity_p}
    With fixed \(Q \succ 0\) and \(Z_i\), define the quadratic form \(V_i(p, y) \coloneqq y^T M_i(p) y\), where \(M_i(p)\) is the left-hand side of either \eqref{eq:ms_stability_decomposable_lmi}, \eqref{eq:mjls_h2_decoupled_gramian_lmi} or \eqref{eq:mjls_h2_decoupled_trace_lmi}.
    \(V_i(p, y)\) is convex in \(p\) for all \(i \in \{1, \ldots, N\}\) if and only if either the relevant matrices from \(\{A^c, B^c, C^c, D^c\}\) are zero matrices or all non-zero eigenvalues of \(L^0\) satisfy \(\lambda_i \geq 2\).
\end{lemma}
\begin{proof}
    We prove the lemma for \eqref{eq:mjls_h2_decoupled_gramian_lmi} as representative of all three inequalities.
    \(V_i(p, y)\) is convex in \(p\) if and only if \cite{Horn2012}
    \begin{equation*}
        \pd[2]{V_i(p, y)}{p} \geq 0 \, \forall y \Leftrightarrow \big(\lambda_i^2 - 2\lambda_i\big)\big(A^{cT} Q A^c + C^{cT} C^c\big) \succeq 0.
    \end{equation*}
    From \(Q \succ 0\), it follows that \(A^{cT} Q A^c \succeq 0\) and \(C^{cT} C^c \succeq 0\) for all \(A^c, C^c\).
    We then distinguish two cases:
    If \(\lambda_i^2 \geq 2\lambda_i\), we are done.
    This condition is satisfied for \(\lambda_i = 0\) and otherwise equivalent to \(\lambda_i \geq 2\) since all eigenvalues of \(L^0\) are non-negative \cite{Mesbahi2010}.
    On the other hand, if \(\lambda_i^2 < 2\lambda_i\), we must have \(A^{cT} Q A^c + C^{cT} C^c \preceq 0\) and thus \(A^{cT} Q A^c = C^{cT} C^c = 0\), which in turn implies \(A^c = 0\) and \(C^c = 0\).
    The proof for \eqref{eq:ms_stability_decomposable_lmi} and \eqref{eq:mjls_h2_decoupled_trace_lmi} follows along the same lines, replacing \((A^c, C^c)\) by \((B^c, D^c)\) for \eqref{eq:mjls_h2_decoupled_trace_lmi} and considering just \(A^c\) for \eqref{eq:ms_stability_decomposable_lmi}.
\end{proof}

For a convex function \(V(x)\), its sublevel set \(\{x : V (x) < 0\}\) is convex as well.
On the interval \([\ubar{p}, \bar{p}]\), this ensures that checking the condition on the boundary is sufficient to verify it is satisfied throughout. Thus, assuming that the conditions of Lemma~\ref{lem:lmi_convexity_p} are fulfilled, the problem is reduced to applying Theorem~\ref{thm:ms_stability_decoupled} or \ref{thm:mjls_h2_decoupled} at \(\ubar{p}\) and \(\bar{p}\) with shared \(Q\) and \(Z_i\).

To give some meaning to the conditions from the lemma, the zero matrix condition implies that there exists no coupling between the agents and is thus irrelevant for the analysis of MAS in practice.
The remaining condition on the eigenvalues of \(L^0\) can be seen as a lower bound on the connectivity of the underlying communication graph \(\mathcal{G}^0\).
In particular, if \(\mathcal{G}^0\) is undirected and connected, it is a lower bound on the Fiedler eigenvalue~\(\lambda_2\), the smallest non-zero eigenvalue of \(L^0\) \cite{Mesbahi2010}.

\subsection{Handling Uncertain Nominal Communication Graphs}\label{sec:scalable_graphs}
In the form stated above, Theorems~\ref{thm:ms_stability_decoupled} and \ref{thm:mjls_h2_decoupled} require complete knowledge of the spectrum of \(L^0\) and thus centralized information.
However, it is possible to utilize another convexity property of the LMIs to relax this restriction.

\begin{lemma}[Convexity in \(\lambda_i\)]\label{lem:lmi_convexity_lambda}
    With fixed \(Q \succ 0\) and \(Z_i\), define the quadratic form \(V_p(\lambda_i, y) \coloneqq y^T M_p(\lambda_i) y\), where \(M_p(\lambda_i)\) is the left-hand side of either \eqref{eq:ms_stability_decomposable_lmi}, \eqref{eq:mjls_h2_decoupled_gramian_lmi} or \eqref{eq:mjls_h2_decoupled_trace_lmi}.
    \(V_p(\lambda_i, y)\) is convex in \(\lambda_i\) for all \(p\).
\end{lemma}
\begin{proof}
    The proof is analogous to the proof of Lemma~\ref{lem:lmi_convexity_p}, outlined exemplarily for \eqref{eq:mjls_h2_decoupled_gramian_lmi}.
    \(Q \succ 0\) guarantees that
    \begin{equation*}
        \pd[2]{V_p(\lambda_i, y)}{\lambda_i} \geq 0 \, \forall y \Leftrightarrow p^2 \big(A^{cT} Q A^c + C^{cT} C^c\big) \succeq 0
    \end{equation*}
    is always satisfied regardless of \(p\).
\end{proof}

Lemma~\ref{lem:lmi_convexity_lambda} implies that knowledge of the boundary of the spectrum of \(L^0\) is sufficient to evaluate Theorems~\ref{thm:ms_stability_decoupled} and \ref{thm:mjls_h2_decoupled}.
An upper bound on \(\lambda_N\) can for example be obtained from the maximum node degree and Cheeger's inequality could be used to bound \(\lambda_2\) \cite[Section 2.4.2]{Mesbahi2010}.
For Theorem~\ref{thm:ms_stability_decoupled}, this adaptation comes without additional conservatism, giving sufficient stability conditions independent of network size.
For Theorem~\ref{thm:mjls_h2_decoupled} on the other hand, one needs to further restrict \(Z_i = Z_j\) for all \(i, j \in \{1 \ldots N\}\), making the upper bound on the \(H_2\)-norm possibly more conservative.

\subsection{Necessary Conditions for the Analysis of MJLS}\label{sec:scalable_necessary}
To evaluate the conservatism introduced by restricting Theorems~\ref{thm:ms_stability_decoupled} and \ref{thm:mjls_h2_decoupled} to a single \(Q\) for all modal constraints, we can compare their results to those obtained from the lossless theorems from Section~\ref{sec:problem_mjls}.
However, this comparison is only tractable for MAS with few agents because of the exponential scaling of the lossless theorems.
Thus, we propose necessary conditions for Theorems~\ref{thm:ms_stability} and \ref{thm:mjls_h2} that can be checked with the same (linear) complexity as the sufficient conditions from Section~\ref{sec:scalable_decomposition}, which enables us to estimate the conservatism for large MAS.

For the analysis, we introduce the \emph{mean system} \(\bar{G}\), which is the LTI system whose system matrices are given by the mean of the MJLS matrices.
For the MJLS~\eqref{eq:mjls_decomposable}, this results in the LTI state-space model
\begin{equation}\label{eq:necessary_lti}
    \bar{G}: \,\left\{\enspace\begin{aligned}
        \bar{x}(k+1) =& \left(I_N \otimes A^d + p L^0 \otimes A^c\right) \bar{x}(k) \\
                     &+ \left(I_N \otimes B^d + p L^0 \otimes B^c\right) \bar{w}(k) \\
        \bar{z}(k)   =& \left(I_N \otimes C^d + p L^0 \otimes C^c\right) \bar{x}(k) \\
                     &+ \left(I_N \otimes D^d + p L^0 \otimes D^c\right) \bar{w}(k).
    \end{aligned}\right.
\end{equation}
The mean system can be seen as advancing the ensemble average state in time, in contrast to the MJLS, which advances one specific realization.
Based on the mean system, we can then state the following result:

\begin{theorem}[LTI Necessary Conditions]\label{thm:necessary_conditions}
    Given the decomposable MJLS~\eqref{eq:mjls_decomposable}, its mean \(\bar{G}\), and any \(\gamma > 0\), the following implications hold:
    \begin{enumerate}[i)]
        \item \(\hat{G}\) is MSS \(\enspace\Rightarrow\enspace\) \(\bar{G}\) is stable
        \item \(\|\hat{G}\|_{H_2} < \gamma \enspace\Rightarrow\enspace \|\bar{G}\|_{H_2} < \gamma\)
    \end{enumerate}
\end{theorem}
\begin{proof}
    See Appendix~\ref{app:theorem_necessary}.
\end{proof}

Theorem~\ref{thm:necessary_conditions} implies that stability of \(\bar{G}\) is necessary for MSS of \(\hat{G}\) and that \(\|\bar{G}\|_{H_2}\) is a lower bound for \(\|\hat{G}\|_{H_2}\).
A similar result can be obtained for general MJLS in that stability in the second moment, i.e. mean-square stability, implies stability in the first moment \cite[Proposition 3.6]{Costa2005}.
Since \(\bar{G}\) is LTI, we may apply the analysis based on modal subsystems proposed by Massioni and Verhaegen in \cite{Massioni2009}, resulting in LMI conditions that scale linearly in the number of agents.

\section{Example: First Order Consensus}\label{sec:example}
\subsection{Setting up the Problem}\label{sec:example_problems}
To demonstrate the scalability of and judge the amount of conservatism in the analysis conditions from Section~\ref{sec:scalable}, let us now finally apply the approach to a numerical example.
The example we chose is the \emph{discrete-time first-order consensus problem}, which can be described as the problem of reaching agreement in a network of linear first-order integrators while each agent is only communicating to a subset of the remaining agents.
The communication between agents is modelled using the graph~\(\mathcal{G}^0\), with agent~\(i\) receiving information according to its in-neighbourhood~\(\mathcal{N}_i^-\).
For each individual agent, the dynamics are then described by \eqref{eq:consensus_agents} with \(x_i(k) \in \real\), while a solution to the consensus problem is given by the \emph{consensus protocol} in \eqref{eq:consensus_protocol} with parameter \(\kappa > 0\) chosen small enough.
\begin{align}
    x_i(k+1) &= x_i(k) + u_i(k) + w_i(k)\label{eq:consensus_agents} \\
    u_i(k)   &= \kappa \sum_{v_j \in \mathcal{N}_i^-} \big(x_j(k) - x_i(k)\big) \label{eq:consensus_protocol}
\end{align}
In contrast to the standard consensus problem, we introduce disturbance inputs \(w_i(k)\) in order to use the \(H_2\)-norm as a performance measure.
By stacking the states and inputs in \(x(k) \coloneqq [x_1(k), \ldots, x_N(k)]^T\) and \(w(k)\), we can write the network dynamics as \(x(k+1) = (I - \kappa L^0) x(k) + w(k)\).
For details on the consensus problem and its solution, see \cite{OlfatiSaber2007}.

For our example, we study the consensus problem with uncertain exchange of information.
Thus, we introduce packet loss using stochastic processes \(\{\alpha_{ij}(k)\}\) adhering to Assumption~\ref{ass:bernoulli} and modify the consensus protocol \eqref{eq:consensus_protocol} to
\begin{equation}\label{eq:consensus_protocol_packet_loss}
    u_i(k) = \kappa \sum_{v_j \in \mathcal{N}_i^-} \alpha_{ij}(k) \big(x_j(k) - x_i(k)\big).
\end{equation}
For the network, this results in \(x(k+1) = (I - \kappa L_\sigma) x(k) + w(k)\), which has the form of the decomposable MJLS~\eqref{eq:mjls_decomposable} with
\begin{align*}
    A^d &=  1,      & B^d &= 1, & C^d &= 1, & D^d &= 0, \\
    A^c &= -\kappa, & B^c &= 0, & C^c &= 0, & D^c &= 0,
\end{align*}
where we are using the full state as performance output, so \(z(k) = x(k)\), and use \(\kappa = 0.1\) in the following.

\begin{figure}
    \centering
    \subfloat[Triangle-shaped graphs \(\mathcal{G}_h^\triangle\)]{
        \centering
        \input{figures/triangle_graphs}
        \label{fig:example_graphs_triangle}
    }
    \hfill
    \subfloat[Circular graphs \(\mathcal{G}_N^\circ\)]{
        \centering
        \input{figures/circular_graphs}
        \label{fig:example_graphs_circular}
    }
    \caption{Graph structures that are used to evaluate the scalability of the analysis conditions.}
    \label{fig:example_graphs}
\end{figure}
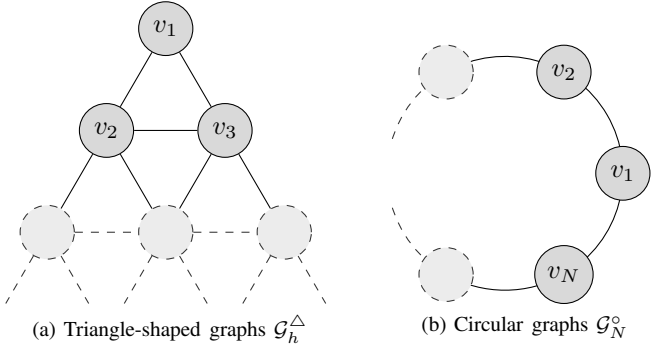

In the numerical example, we are using the two graph structures shown in Fig.~\ref{fig:example_graphs} to interconnect the agents.
The family of circular graphs \(\mathcal{G}_N^\circ\) shown in Fig.~\ref{fig:example_graphs_circular} has twice as many edges as vertices (note that \(e_{ij}\) is different from \(e_{ji}\)), which makes it suitable for testing the scalability of Theorem~\ref{thm:mjls_h2} because applying the theorem to networks even with edge counts in the low double digits is challenging.
On the other hand, they suffer from poor connectivity for networks with many agents.
Therefore, we will be using the triangle shaped graphs~\(\mathcal{G}_h^\triangle\) from Fig.~\ref{fig:example_graphs_triangle} for larger networks.
If the number of vertices in the last row is denoted by \(h\), then \(\mathcal{G}_h^\triangle\) has \(N = \frac{h}{2} (h+1)\) vertices and \(|\mathcal{E}_h^\triangle| = 3h (h-1)\) edges.
Note that both graph structures are undirected and thus satisfy Assumption~\ref{ass:undirected_graph}.

An issue with analysing the system that describes the solution to the consensus problem is its inherent marginal stability.
In the decomposable systems framework without packet loss from \cite{Massioni2009}, a convenient approach to resolve this issue is to neglect the modal subsystem that corresponds to the 0 eigenvalue of the Laplacian \cite{Fax2004}.
For the calculation of the \(H_2\)-norm, this means that \((\Pi \otimes I_{n_z}) \hat{G}\) is analysed instead of \(\hat{G}\), where \(\Pi \coloneqq I_N - \frac 1N \mathbf{1}_N \mathbf{1}_N^T\) is the orthogonal projection onto the disagreement space.
The same approach can be applied to the decoupled analysis conditions from Theorems~\ref{thm:ms_stability_decoupled} and \ref{thm:mjls_h2_decoupled} as well as the necessary LTI conditions in Theorem~\ref{thm:necessary_conditions} and -- in adapted form -- the coupled LMIs in Theorems~\ref{thm:ms_stability} and \ref{thm:mjls_h2}.

Consider again a transformation \(U\) with \(U^T U = I_N\) such that \(U^T L^0 U\) is diagonal, which does exist under Assumption~\ref{ass:undirected_graph}.
Note that because \(L^0\) has zero row and column sum, \(U\) can be chosen as \(U = [\mathbf{1}_N / \sqrt{N} \; \tilde{U}]\) with \(\Pi U = [0 \; \tilde{U}]\).
We can then apply \(U\) as state and signal transformation to \((\Pi \otimes I_{n_z}) \hat{G}\), giving \(\tilde{x}(k) \coloneqq (U^T \otimes I_{n_x}) x(k)\), \(\tilde{w}(k) \coloneqq (U^T \otimes I_{n_w}) w(k)\) and \(\tilde{z}(k) \coloneqq ((U^T \Pi) \otimes I_{n_z}) z(k)\), resulting in
\begin{subequations}\label{eq:mjls_transformed}\begin{align}\begin{split}
        \tilde{x}(k+1) =& \bbma A^d & \tilde{l}_i \otimes A^c \\ 0 & I_{N-1} \otimes A^d + \tilde{L}_i \otimes A^c \ebma \tilde{x}(k) \\
                       &+ \bbma B^d & \tilde{l}_i \otimes B^c \\ 0 & I_{N-1} \otimes B^d + \tilde{L}_i \otimes B^c \ebma \tilde{w}(k),
    \end{split}\\
    \begin{split}
        \tilde{z}(k)   =& \bbma 0 & 0 \\ 0 & I_{N-1} \otimes C^d + \tilde{L}_i \otimes C^c \ebma \tilde{x}(k) \\
                       &+ \bbma 0 & 0 \\ 0 & I_{N-1} \otimes D^d + \tilde{L}_i \otimes D^c \ebma \tilde{w}(k),
\end{split}\end{align}\end{subequations}
where \(\tilde{l}_i \coloneqq \mathbf{1}_N^T L_i \tilde{U} / \sqrt{N}\) and \(\tilde{L}_i \coloneqq \tilde{U}^T L_i \tilde{U}\).
It was shown in \cite{Massioni2009} that the \(H_2\)-norm is invariant under this kind of orthogonal transformation of input and output.

In the transformed system~\eqref{eq:mjls_transformed}, it is apparent that the centre of gravity, which \(\tilde{x}_k\) contains in its first \(n_x\) entries, does not affect the remaining states, since the bottom left block of every system matrix is 0.
The converse does however only hold if \emph{all} \(\mathcal{G}_i\) are balanced, since this implies that \(\tilde{l}_i = 0\) for all \(i \in \mathcal{K}\).
For stability analysis using Theorem~\ref{thm:ms_stability}, this one-way coupling may be ignored, as stability of the remaining system would imply the centre of gravity stays finite as long as the decoupled part \(A^d\) is at least marginally stable.
On the other hand, when calculating the \(H_2\)-norm of the system, we can take advantage of the fact that the full first columns are zero for the \(C\) and \(D\) matrices of the transformed system.
This implies that even though the centre of gravity is affected by the remaining system, this is not apparent in the output \(\tilde{z}(k)\) and thus does not increase the \(H_2\)-norm.
The desired \(H_2\)-norm can therefore be obtained by only considering the bottom right block of the transformed system.

\subsection{Numerical results}\label{sec:example_results}
To analyse the scalability and conservatism of the approaches described in this paper, we implemented the LMI conditions from Theorems~\ref{thm:mjls_h2}, \ref{thm:mjls_h2_decoupled}, and \ref{thm:necessary_conditions} in \textsc{Matlab} using the \textsc{Yalmip} \cite{Loefberg2004} toolbox.
All three conditions are affine in \(\gamma^2\) such that we can directly minimize \(\gamma^2\) -- and thus \(\gamma\) -- subject to either of the LMIs.
The minimum \(\gamma\) obtainable by each of the conditions will be plotted below as the respective \(H_2\)-performance.
All source code is available at \cite{SourceCode}.

\begin{figure}
    \centering
    \subfloat[Network with communication graph \(\mathcal{G}_3^\triangle\)]{
        \centering
        \input{figures/h2_triangle_small}
        \label{fig:example_conservatism_small}
    }
    \\
    \subfloat[Network with communication graph \(\mathcal{G}_{50}^\triangle\)]{
        \centering
        \begin{tikzpicture}
            \pgfplotstableread[col sep=comma]{figures/data/conservatism_large.csv} \datatable

            \begin{semilogyaxis}[
                width=0.9\columnwidth,
                xlabel=Transmission probability $p$,
                ylabel=$H_2$-performance,
                ymax=2e4,
            ]
                \addplot table[x=p, y=large_decom] {\datatable};
                \addplot table[x=p, y=large_mean] {\datatable};

                \legend{Decomposed, Mean};
            \end{semilogyaxis}
        \end{tikzpicture}
        \label{fig:example_conservatism_large}
    }
    \caption{\(H_2\)-performance of the MJLS at different transmission probabilities~\(p\) with each of the three analysis conditions.}
    \label{fig:example_conservatism}
\end{figure}
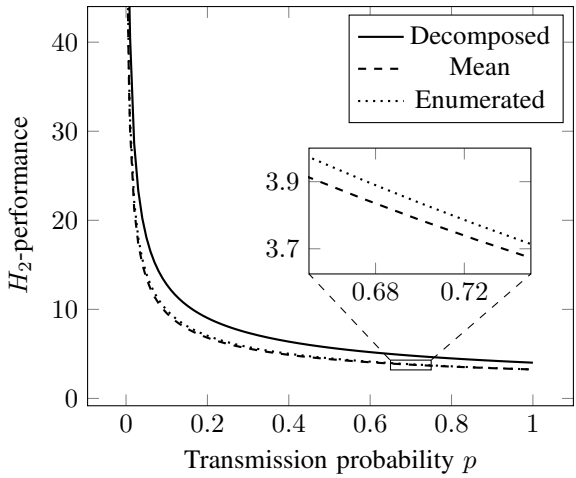
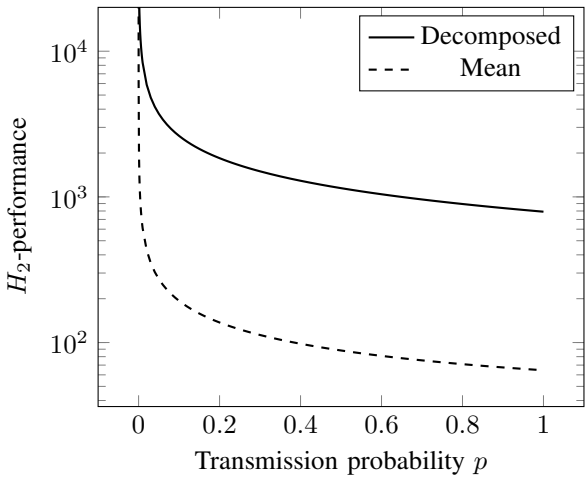

Let us first evaluate how the \(H_2\)-performance changes with the transmission probability~\(p\) for each of the three conditions.
For two test graphs, \(\mathcal{G}_3^\triangle\) and \(\mathcal{G}_{50}^\triangle\), we perform a sweep over \(p\), which is shown in Fig.~\ref{fig:example_conservatism}.
In the figure, the \enquote{\emph{decomposed}} graph refers to the best upper bound on the \(H_2\)-norm that can be obtained from Theorem~\ref{thm:mjls_h2_decoupled}, \enquote{\emph{mean}} is the lower bound based on Theorem~\ref{thm:necessary_conditions} and \enquote{\emph{enumerated}} corresponds to the original analysis condition in Theorem~\ref{thm:mjls_h2} and thus shows the true \(H_2\)-norm of the system.
Theorem~\ref{thm:mjls_h2} is only applied to the small MAS in Fig.~\ref{fig:example_conservatism_small}, since \(\mathcal{G}_{50}^\triangle\) has \(m = 2^{|\mathcal{E}_{50}^\triangle|} = 2^{7350}\) modes, which are intractable to enumerate.

As expected, the performance figures obtained from the decomposed analysis results in Theorems~\ref{thm:mjls_h2_decoupled} and \ref{thm:necessary_conditions} do not match the \(H_2\)-norm of the system but over- and underestimate it, respectively.
Furthermore, the gap between the upper and lower bound is significantly increased for the larger MAS.
However, while the mean system recovers the exact norm for \(p = 1\) because it coincides with the MJLS, the upper bound is conservative for all transmission probabilities.

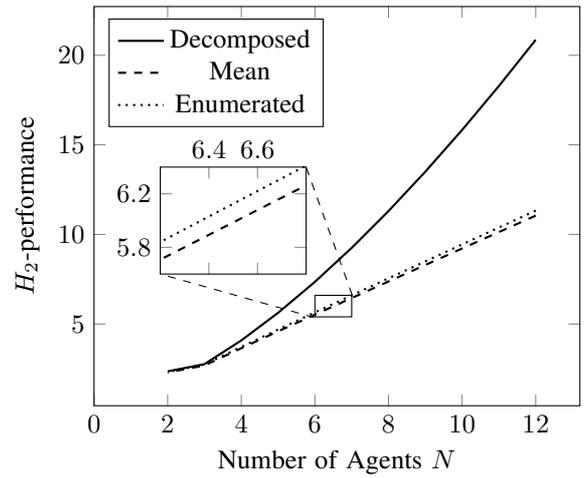
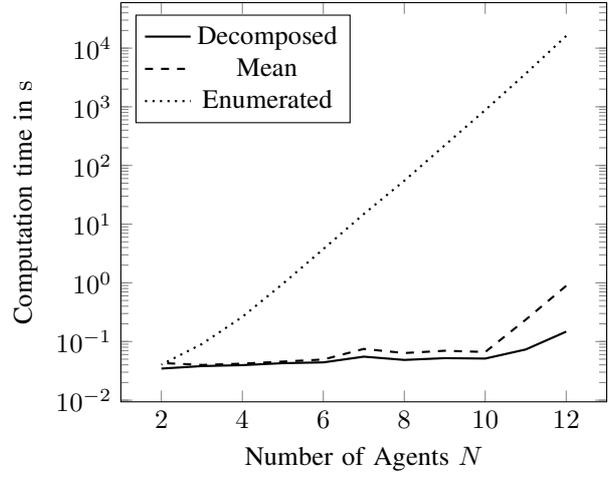
\begin{figure}
    \pgfplotstableread[col sep=comma]{figures/data/scaling_small.csv} \datatable

    \centering
    \subfloat[\(H_2\)-performance of the MJLS at \(p = 0.5\)]{
        \centering
        \input{figures/h2_circle_small}
        \label{fig:example_scaling_small_norm}
    }
    \\
    \subfloat[Computation Time]{
        \centering
        \begin{tikzpicture}
            \begin{semilogyaxis}[
                width=0.9\columnwidth,
                xlabel=Number of Agents $N$,
                ylabel=Computation time in \si{\second},
                legend pos=north west,
            ]
                \addplot table[x=n, y=time_decom] {\datatable};
                \addplot table[x=n, y=time_mean] {\datatable};
                \addplot table[x=n, y=time_enum] {\datatable};

                \legend{Decomposed, Mean, Enumerated};
            \end{semilogyaxis}
        \end{tikzpicture}
        \label{fig:example_scaling_small_time}
    }
    \caption{Scalability of the three analysis conditions for small networks with graphs \(\mathcal{G}_N^\circ\).}
    \label{fig:example_scaling_small}
\end{figure}

In a second step, we compare how the analysis conditions from Theorems~\ref{thm:mjls_h2} and \ref{thm:mjls_h2_decoupled} scale in terms of computational speed and conservatism of the calculated \(H_2\)-norm.
We start by analysing the MAS with the circular graphs \(\mathcal{G}_N^\circ\) for \(N\) between 2 and 12.
Because the number of edges is relatively small for these graphs, we can apply all three conditions.
For each \(N\), we calculate the \(H_2\)-performance with \(p = 0.5\), which was chosen since it is the transmission probability with the largest variance.
The results are shown in Fig.~\ref{fig:example_scaling_small}.

As observed before, the conservatism of Theorem~\ref{thm:mjls_h2_decoupled} grows with increasing agent count.
The lower bound obtained from the mean system is close to the exact norm regardless of the agent count.
Concerning the computational speed, it is apparent that the analysis conditions from Theorem~\ref{thm:mjls_h2} show an exponential growth in complexity such that the problem will quickly become intractable even for networks of moderate size.
On the other hand, the conditions from Theorem~\ref{thm:mjls_h2_decoupled} show no substantial increase in computation time.

\begin{figure}
    \pgfplotstableread[col sep=comma]{figures/data/scaling_large.csv} \datatable

    \centering
    \subfloat[\(H_2\)-performance of the MJLS at \(p = 0.5\)]{
        \centering
        \begin{tikzpicture}
            \begin{loglogaxis}[
                width=0.9\columnwidth,
                xlabel=Number of Agents $N$,
                ylabel=$H_2$-performance,
                legend pos=north west,
            ]
                \addplot table[x=N, y=norm_decom] {\datatable};
                \addplot table[x=N, y=norm_mean] {\datatable};

                \legend{Decomposed, Mean};
            \end{loglogaxis}
        \end{tikzpicture}
        \label{fig:example_scaling_large_norm}
    }
    \\
    \subfloat[Computation Time]{
        \centering
        \begin{tikzpicture}
            \begin{loglogaxis}[
                width=0.9\columnwidth,
                xlabel=Number of Agents $N$,
                ylabel=Computation time in \si{\second},
                legend pos=north west,
                log y ticks with fixed point,
            ]
                \addplot table[x=N, y=time_decom] {\datatable};
                \addplot table[x=N, y=time_mean] {\datatable};

                \legend{Decomposed, Mean};
            \end{loglogaxis}
        \end{tikzpicture}
        \label{fig:example_scaling_large_time}
    }
    \caption{Scalability of the two fast analysis conditions for large networks with graphs \(\mathcal{G}_h^\triangle\).}
    \label{fig:example_scaling_large}
\end{figure}
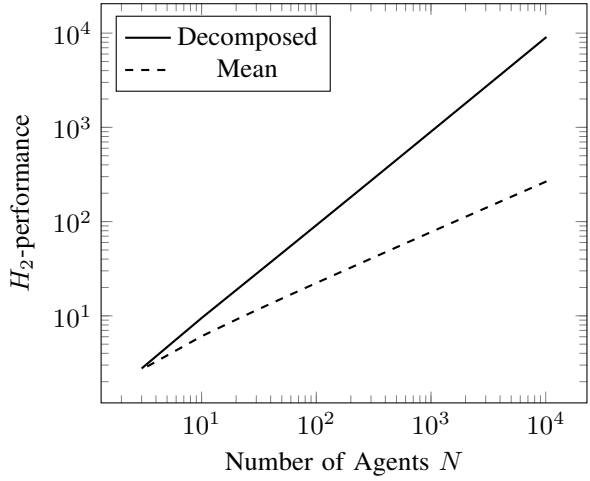
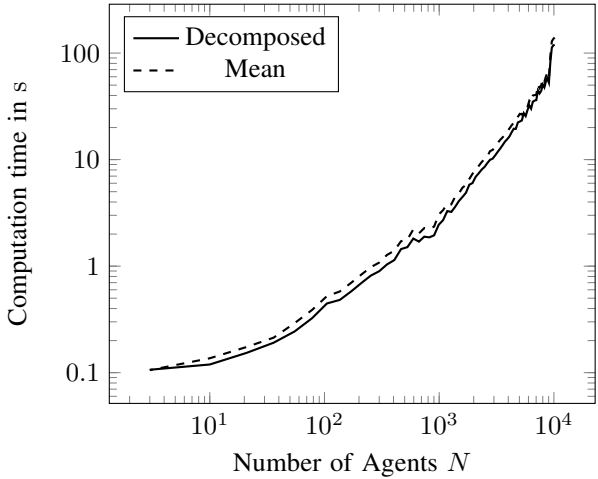

Finally, for the last test we apply the decomposed analysis conditions to MAS based on the triangle-shaped graphs \(\mathcal{G}_2^\triangle\) to \(\mathcal{G}_{142}^\triangle\), ranging from a network with three agents to one with 10011.
Again, we calculate bounds on the \(H_2\)-norm for \(p = 0.5\) with each of the graphs.
The corresponding performance and computation time curves are displayed in Fig.~\ref{fig:example_scaling_large}.

While the exact \(H_2\)-norm for these systems is unknown, the gap between the upper and lower bound on the performance widens with increasing number of agents, up to one and a half orders of magnitude for the largest MAS.
Still, the decomposed conditions allow to calculate an upper bound on the \(H_2\)-performance for networks of that size at all, in contrast to the original conditions from Theorem~\ref{thm:mjls_h2}, which is intractable to validate for systems that large.
In terms of computational speed, the linear scaling of the decomposed conditions is confirmed by Fig.~\ref{fig:example_scaling_large_time} for networks with more than 40 agents.
For smaller MAS the constant cost of setting up the optimization problem shadows the cost of solving the LMIs, resulting in diminishing returns when the number of agents is decreased below that threshold.

\section{Conclusions and Future Work}\label{sec:conclusions}
This paper proposes an extension of the decomposable systems framework to stochastic jump linear systems in order to analyse the effect of Bernoulli distributed packet loss with uniform packet loss probability on multi-agent systems.
Based on analytic expressions for the expected Laplacians, sufficient analysis conditions for mean-square stability and bounds on the \(H_2\)-norm that scale linearly with the number of agents were derived.
Finally, it was demonstrated that the proposed conditions are applicable to very large networks but that their conservatism increases with the size of the network.

In future work, it will be investigated if the restriction to identical matrix variables in the presented analysis conditions can be removed without losing sufficiency, possibly leading to lossless complexity reduction similar to the LTI case.
An instrumental step would be to extend the result on simultaneous diagonalizability to the more general weighted expectation.
Furthermore, current research is aiming at how the restrictive assumption of Bernoulli distributed loss with uniform probability can be relaxed.

%% file: figures/triangle_graphs.tex
\begin{tikzpicture}[scale=0.9]
    \pgfmathsetmacro{\ra}{cos(pi/6)}
    \pgfmathsetmacro{\rb}{2*\ra}

	\draw ( 30:\ra) node[agent](a3){$v_3$};
	\draw ( 90:\rb) node[agent](a1){$v_1$};
	\draw (150:\ra) node[agent](a2){$v_2$};
	\draw (210:\rb) node[hidden agent](a4){\phantom{$v_1$}};
	\draw (270:\ra) node[hidden agent](a5){\phantom{$v_1$}};
	\draw (330:\rb) node[hidden agent](a6){\phantom{$v_1$}};
	
	\draw[link] (a1) -- (a2);
	\draw[link] (a1) -- (a3);
	\draw[link] (a2) -- (a3);
	\draw[link] (a2) -- (a4);
	\draw[link] (a2) -- (a5);
	\draw[link] (a3) -- (a5);
	\draw[link] (a3) -- (a6);
	\draw[hidden link] (a4) -- (a5);
	\draw[hidden link] (a5) -- (a6);
	
	\draw[hidden link] (a4) -- ++($0.7*(a4)-0.7*(a2)$);
	\draw[hidden link] (a5) -- ++($0.7*(a4)-0.7*(a2)$);
	\draw[hidden link] (a6) -- ++($0.7*(a4)-0.7*(a2)$);
	\draw[hidden link] (a4) -- ++($0.7*(a6)-0.7*(a3)$);
	\draw[hidden link] (a5) -- ++($0.7*(a6)-0.7*(a3)$);
	\draw[hidden link] (a6) -- ++($0.7*(a6)-0.7*(a3)$);
\end{tikzpicture}

%% file: figures/circular_graphs.tex
\begin{tikzpicture}[scale=1.55]
    \draw[link] (120:1) arc (120:-120:1);
    \draw[hidden link] (120:1) arc (120:165:1);
    \draw[hidden link] (240:1) arc (240:195:1);

	\draw (  0:1) node[agent](a1){$v_1$};
	\draw ( 60:1) node[agent](a2){$v_2$};
	\draw (120:1) node[hidden agent](a3){\phantom{$v_1$}};
	\draw (240:1) node[hidden agent](a4){\phantom{$v_1$}};
	\draw (300:1) node[agent, inner sep=0.3em](a5){$v_N$};
\end{tikzpicture}

%% file: figures/h2_triangle_small.tex
\begin{tikzpicture}
    \pgfplotstableread[col sep=comma]{figures/data/conservatism_small.csv} \datatable
    
    \begin{axis}[
        width=0.9\columnwidth,
        xlabel=Transmission probability $p$,
        ylabel=$H_2$-performance,
        ymax=44,
        restrict y to domain=0:100,
    ]
        \node[draw, fit={(0.65,3.2)(0.75,4.3)}, inner sep=0](magrect){};
    
        \addplot table[x=p, y=small_decom] {\datatable};
        \addplot table[x=p, y=small_mean] {\datatable};
        \addplot table[x=p, y=small_enum] {\datatable};
        
        \legend{Decomposed, Mean, Enumerated};
        
        \coordinate (insetSW) at (axis cs:0.45,14);
    \end{axis}
    \begin{axis}[
        name=inset,
        at={(insetSW)},
        height=3.25cm,
        width=4.5cm,
        xmin=0.65,
        xmax=0.75,
        ymin=3.625,
        ymax=4,
        xtick={0.68,0.72},
        ytick={3.7,3.9},
        restrict y to domain=0:100,
    ]
        \pgfplotsset{cycle list shift=1};
        \addplot table[x=p, y=small_mean] {\datatable};
        \addplot table[x=p, y=small_enum] {\datatable};
    \end{axis} 
    
    \draw[dashed] (magrect.north west) -- (inset.south west);
    \draw[dashed] (magrect.north east) -- (inset.south east);
\end{tikzpicture}

%% file: figures/h2_circle_small.tex
\begin{tikzpicture}
    \begin{axis}[
        width=0.9\columnwidth,
        xlabel=Number of Agents $N$,
        ylabel=$H_2$-performance,
        legend pos=north west,
        xmin=0,
    ]
        \node[draw, fit={(6,5.4)(7,6.6)}, inner sep=0](magrect){};
    
        \addplot table[x=n, y=norm_decom] {\datatable};
        \addplot table[x=n, y=norm_mean] {\datatable};
        \addplot table[x=n, y=norm_enum] {\datatable};
        
        \legend{Decomposed, Mean, Enumerated};
        
        \coordinate (insetSW) at (axis cs:1.8,7.8);
    \end{axis}
    \begin{axis}[
        name=inset,
        at={(insetSW)},
        height=3cm,
        width=3.5cm,
        xmin=6.2,
        xmax=6.8,
        ymin=5.6,
        ymax=6.4,
        xtick={6.4,6.6},
        ytick={5.8,6.2},
        xticklabel pos=top,
    ]
        \pgfplotsset{cycle list shift=1};
        \addplot table[x=n, y=norm_mean] {\datatable};
        \addplot table[x=n, y=norm_enum] {\datatable};
    \end{axis} 
    
    \draw[dashed] (magrect.south west) -- (inset.south west);
    \draw[dashed] (magrect.north east) -- (inset.north east);
\end{tikzpicture}

%% file: doc/appendix.tex
\section{Proof of Lemma~\ref{lem:expected_laplacians}}\label{app:lemma_expectation}
\begin{proof}
    The expectations are calculated element-wise.
    Thus, for \(\expect[L_\sigma]\) we get \(\expect[l_{ij}^\sigma] = -p\) if \(v_j \in \mathcal{N}_i^-\) and \(\expect[l_{ij}^\sigma] = 0\) otherwise for the off-diagonal entries.
    On the diagonal, we have \(\expect[l_{ii}^\sigma] = p d_i^-\).
    Together, this is equal to \(p L(\mathcal{G}^0)\).

    On the other hand, for the expectation of \(L_\sigma^T L_\sigma\), calculate the entries of \(L(\mathcal{G}^0)^T L(\mathcal{G}^0)\) first.
    We get
    \begin{align}
        l_i^{0T} l_i^0 &= (d_i^-)^2 + d_i^+, \label{eq:laplacian_nominal_entries_diag}\\
        l_i^{0T} l_j^0 &= \big|\mathcal{N}_i^+ \cap \mathcal{N}_j^+\big| - d_i^- \mathbb{I}_{\mathcal{N}_i^-}(v_j) - d_j^- \mathbb{I}_{\mathcal{N}_j^-}(v_i) \label{eq:laplacian_nominal_entries_offdiag}
    \end{align}
    for the diagonal and off-diagonal entries respectively, where \(l_i^0\) and \(l_j^0\) are the \(i\)th and \(j\)th column of \(L^0\).
    Notice that \(v_i \in \mathcal{N}_s^- \Leftrightarrow v_s \in \mathcal{N}_i^+\) and \(\{v_i, v_j\} \subseteq \mathcal{N}_s^- \Leftrightarrow v_s \in \mathcal{N}_i^+ \cap \mathcal{N}_j^+\).
    Then, define \(\beta_{ij}\) as the elements of \(\expect[L_\sigma^T L_\sigma]\) and see that
    \begin{equation}\label{eq:laplacian_expectation_entries}
        \beta_{ij} = \expect\left[l_i^{\sigma T} l_j^\sigma\right] = \sum_{s = 1}^N \expect[l_{si}^\sigma l_{sj}^\sigma] \eqqcolon \sum_{s = 1}^N \beta_{ij}^s.
    \end{equation}
    To calculate their values, recall the definition of \(l_{ij}^\sigma\) in \eqref{eq:random_laplacian} and distinguish the following five cases:
    \begin{equation*}\beta_{ij}^s = \begin{cases}
            p \mathbb{I}_{\mathcal{N}_i^+}(v_s)                                   & if \(i = j \neq s\), \\
            (d_i^-)^2 p^2 + d_i^- p(1-p) \vphantom{\mathbb{I}_{\mathcal{N}_i^+}}  & if \(i = j = s\), \\
            p^2 \mathbb{I}_{\mathcal{N}_i^+ \cap \mathcal{N}_j^+}(v_s)            & if \(s \neq i \neq j \neq s\), \\
            -\big(d_i^- p^2 + p(1-p)\big) \mathbb{I}_{\mathcal{N}_i^-}(v_j)       & if \(s = i \neq j\), \\
            -\big(d_j^- p^2 + p(1-p)\big) \mathbb{I}_{\mathcal{N}_j^-}(v_i)       & if \(i \neq j = s\).
    \end{cases}\end{equation*}
    With all five cases covered, sum up the results according to \eqref{eq:laplacian_expectation_entries} to obtain \(\beta_{ij}\).
    On the main diagonal, we have
    \begin{equation}\label{eq:laplacian_expectation_entries_diag}
        \beta_{ii} = p^2 \big((d_i^-)^2 + d_i^+\big) + p(1-p) \big(d_i^- + d_i^+\big),
    \end{equation}
    while for the off-diagonal entries use \(\mathbb{I}_{\mathcal{N}_j^-}(v_i) = \mathbb{I}_{\mathcal{N}_i^+}(v_j)\) to arrive at
    \begin{equation}\label{eq:laplacian_expectation_entries_offdiag}\begin{aligned}
            \beta_{ij} =& p^2 \Big(\big|\mathcal{N}_i^+ \cap \mathcal{N}_j^+\big| - d_i^- \mathbb{I}_{\mathcal{N}_i^-}(v_j) - d_j^- \mathbb{I}_{\mathcal{N}_j^-}(v_i) \Big) \\
            &- p(1-p) \Big(\mathbb{I}_{\mathcal{N}_i^-}(v_j) + \mathbb{I}_{\mathcal{N}_i^+}(v_j)\Big).
    \end{aligned}\end{equation}
    Notice that \eqref{eq:laplacian_expectation_entries_diag} and \eqref{eq:laplacian_expectation_entries_offdiag} contain \(p^2\) multiplied by \eqref{eq:laplacian_nominal_entries_diag} and \eqref{eq:laplacian_nominal_entries_offdiag} respectively, resulting in the first term in the lemma.
    Finally, use the fact that in the transposed graph \(\mathcal{G}^{0T}\) the in- and out-neighbourhoods are exchanged compared to the original graph \(\mathcal{G}^0\) to see that the remaining terms correspond to the second part of the equation.
\end{proof}
\begin{remark}
    It is possible to exclude opposing links from the independence clause in Assumption~\ref{ass:bernoulli} because there are no products between \(\alpha_{ij}\) and \(\alpha_{ji}\) for any pair \((i,j)\) in the calculations leading up to \(\beta_{ij}^s\).
\end{remark}

\section{Proof of Theorem~\ref{thm:mjls_h2_decoupled}}\label{app:theorem_decomposed_h2}
\begin{proof}
    First, notice that \eqref{eq:mjls_h2_gramian_lmi} is equivalent to
    \begin{equation}\label{eq:mjls_h2_gramian_lmi_expectation}
        \expect\left[A_\sigma^T Q A_\sigma + C_\sigma^T C_\sigma\right] - Q \prec 0.
    \end{equation}
    Then, imposing \(Q = I_N \otimes \tilde{Q}\), apply the same steps as in the proof to Corollary~\ref{cor:ms_stability_decomposable} to arrive at
    \begin{align*}
        &I_N \otimes \big(A^{dT} \tilde{Q} A^d + C^{dT} C^d - \tilde{Q}\big) \\
        &+ \expect\big[L_\sigma^T\big] \otimes \big(A^{cT} \tilde{Q} A^d + C^{cT} C^d\big) \\
        &+ \expect\big[L_\sigma\big] \otimes \big(A^{dT} \tilde{Q} A^c + C^{dT} C^c\big) \\
        &+ \expect\big[L_\sigma^T L_\sigma\big] \otimes \big(A^{cT} \tilde{Q} A^c + C^{cT} C^c\big) \prec 0.
    \end{align*}
    Following the proof of Theorem~\ref{thm:ms_stability_decoupled}, we utilize Lemmas~\ref{lem:expected_laplacians} and \ref{lem:simultaneous_diagonalization} to apply a congruence transformation, resulting in
    \begin{align*}
        &I_N \otimes \big(A^{dT} \tilde{Q} A^d + C^{dT} C^d - \tilde{Q}\big) \\
        &+ \big(p^2\Lambda^2 + 2p(1-p) \Lambda\big) \otimes \big(A^{cT} \tilde{Q} A^c + C^{cT} C^c\big) \\
        &+ p\Lambda \otimes \big(A^{dT} \tilde{Q} A^c + C^{dT} C^c + A^{cT} \tilde{Q} A^d + C^{cT} C^d\big) \prec 0.
    \end{align*}
    Since every component is block diagonal, this is equivalent to \eqref{eq:mjls_h2_decoupled_gramian_lmi}.
    Moreover, apply the same steps to \eqref{eq:mjls_h2_trace_lmi} without imposing additional constraints on \(Z\), leading to
    \begin{align*}
        &I_N \otimes \big(B^{dT} \tilde{Q} B^d + D^{dT} D^d\big) - \tilde{Z}\\
        &+ \big(p^2\Lambda^2 + 2p(1-p) \Lambda\big) \otimes \big(B^{cT} \tilde{Q} B^c + D^{cT} D^c\big) \\
        &+ p\Lambda \otimes \big(B^{dT} \tilde{Q} B^c + D^{dT} D^c + B^{cT} \tilde{Q} B^d + D^{cT} D^d\big) \prec 0,
    \end{align*}
    with \(\tilde{Z} \coloneqq (U^T \otimes I_{n_w}) Z (U \otimes I_{n_w})\).
    Neglecting \(\tilde{Z}\), this LMI is block diagonal.
    The Schur complement implies that the diagonal blocks of any negative definite matrix must be negative definite, therefore we can -- without loss of generality -- assume \(\tilde{Z}\) is block diagonal with \(\tilde{Z}_i\) on the diagonal, such that the LMI becomes equivalent to \eqref{eq:mjls_h2_decoupled_trace_lmi}.
    Finally, \(\trace(Z) = \sum_{i = 1}^N \trace\big(\tilde{Z}_i\big)\) and we rename \(\tilde{Q} \to Q\) and \(\tilde{Z}_i \to Z_i\).

    Furthermore, because all \(\lambda_i\) are non-negative \cite{Mesbahi2010} and \(p\) is in the interval \([0,1]\), \(\bar{C}_i^T \bar{C}_i + 2p(1-p)\lambda_i C^{cT}C^c \succeq 0\) and thus any \(Q \succ 0\) satisfying \eqref{eq:mjls_h2_decoupled_gramian_lmi} also satisfies \eqref{eq:ms_stability_decoupled_lmi}.
    This implies mean-square stability for \(\hat{G}\) according to Theorem~\ref{thm:ms_stability_decoupled} and completes the proof by application of Theorem~\ref{thm:mjls_h2}.
\end{proof}

\section{Proof of Theorem~\ref{thm:necessary_conditions}}\label{app:theorem_necessary}
\begin{proof}
    Both implications in the theorem are based on showing that the LMIs in Theorems~\ref{thm:ms_stability} and \ref{thm:mjls_h2} imply their LTI counterparts for \(\bar{G}\).
    As for Lemma~\ref{lem:lmi_convexity_p}, we demonstrate the argument for \eqref{eq:mjls_h2_gramian_lmi} only, since it can be applied to \eqref{eq:ms_stability_lmi} and \eqref{eq:mjls_h2_trace_lmi} analogously.

    The main idea of the proof is to exploit a definiteness property of a variance like term for matrix-valued random variables.
    For any matrix-valued random variable \(X\), we have
    \begin{align*}
        0 &\preceq \expect\big[(X - \expect[X])^T (X - \expect[X])\big] \\
        &= \expect\big[X^T X - X^T \expect[X] - \expect[X]^T X + \expect[X]^T \expect[X]\big] \\
        &= \expect\big[X^T X\big] - \expect[X]^T \expect[X].
    \end{align*}
    To apply this result, we restate the sum in \eqref{eq:mjls_h2_gramian_lmi} as
    \begin{equation*}
        \sum_{j \in \mathcal{K}} t_j \left(A_j^T Q A_j + C_j^T C_j\right)
        = \expect\big[\tilde{A}_\sigma^T \tilde{A}_\sigma\big] + \expect\big[C_\sigma^T C_\sigma\big]
    \end{equation*}
    where \(\tilde{A}_\sigma \coloneqq Q^{\frac 12} A_\sigma\).
    We thus obtain
    \begin{equation}\label{eq:necessary_lti_expectation}\begin{aligned}
            \expect\big[\tilde{A}_\sigma\big]^T \expect\big[\tilde{A}_\sigma\big] + \expect\big[C_\sigma\big]^T \expect\big[C_\sigma\big] - Q &=\\
            \expect\big[A_\sigma\big]^T Q \expect\big[A_\sigma\big] + \expect\big[C_\sigma\big]^T \expect\big[C_\sigma\big] - Q &\prec 0
    \end{aligned}\end{equation}
    as a necessary condition for \eqref{eq:mjls_h2_gramian_lmi}.
    By Lemma~\ref{lem:expected_laplacians}, we have \(\expect[A_\sigma] = I_N \otimes A^d + p L^0 \otimes A^c\), which corresponds to \(\bar{G}\), and similarly for \(\expect[C_\sigma]\).
    Thus, \eqref{eq:necessary_lti_expectation} is equivalent to the first LMI required to calculate the \(H_2\)-norm of \(\bar{G}\) \cite{Caverly2019}.
    The same procedure can be applied to \eqref{eq:ms_stability_lmi} and \eqref{eq:mjls_h2_trace_lmi}.
\end{proof}

%% file: arxiv.bbl
\begin{thebibliography}{10}
\providecommand{\url}[1]{#1}
\csname url@samestyle\endcsname
\providecommand{\newblock}{\relax}
\providecommand{\bibinfo}[2]{#2}
\providecommand{\BIBentrySTDinterwordspacing}{\spaceskip=0pt\relax}
\providecommand{\BIBentryALTinterwordstretchfactor}{4}
\providecommand{\BIBentryALTinterwordspacing}{\spaceskip=\fontdimen2\font plus
\BIBentryALTinterwordstretchfactor\fontdimen3\font minus
  \fontdimen4\font\relax}
\providecommand{\BIBforeignlanguage}[2]{{%
\expandafter\ifx\csname l@#1\endcsname\relax
\typeout{** WARNING: IEEEtran.bst: No hyphenation pattern has been}%
\typeout{** loaded for the language `#1'. Using the pattern for}%
\typeout{** the default language instead.}%
\else
\language=\csname l@#1\endcsname
\fi
#2}}
\providecommand{\BIBdecl}{\relax}
\BIBdecl

\bibitem{Massioni2009}
P.~Massioni and M.~Verhaegen, ``Distributed control for identical dynamically
  coupled systems: A decomposition approach,'' \emph{{IEEE} Transactions on
  Automatic Control}, vol.~54, no.~1, pp. 124--135, Jan. 2009.

\bibitem{Mesbahi2010}
M.~Mesbahi and M.~Egerstedt, \emph{Graph Theoretic Methods in Multiagent
  Networks}.\hskip 1em plus 0.5em minus 0.4em\relax Princeton University Press,
  Aug. 2010.

\bibitem{Fax2004}
J.~A. Fax and R.~M. Murray, ``Information flow and cooperative control of
  vehicle formations,'' \emph{{IEEE} Transactions on Automatic Control},
  vol.~49, no.~9, pp. 1465--1476, Sep. 2004.

\bibitem{Hoffmann2013}
C.~Hoffmann, A.~Eichler, and H.~Werner, ``Distributed control of linear
  parameter-varying decomposable systems,'' in \emph{American Control
  Conference}.\hskip 1em plus 0.5em minus 0.4em\relax {IEEE}, Jun. 2013.

\bibitem{Eichler2013}
A.~Eichler, C.~Hoffmann, and H.~Werner, ``Robust stability analysis of
  interconnected systems with uncertain time-varying time delays via {IQCs},''
  in \emph{52nd {IEEE} Conference on Decision and Control}.\hskip 1em plus
  0.5em minus 0.4em\relax {IEEE}, Dec. 2013.

\bibitem{Ma2020}
J.~Ma, X.~Yu, and W.~Lan, ``Distributed consensus of linear multi-agent systems
  with nonidentical random packet loss,'' in \emph{59th {IEEE} Conference on
  Decision and Control}.\hskip 1em plus 0.5em minus 0.4em\relax {IEEE}, Dec.
  2020.

\bibitem{Zhang2017}
W.~Zhang, Y.~Tang, T.~Huang, and J.~Kurths, ``Sampled-data consensus of linear
  multi-agent systems with packet losses,'' \emph{{IEEE} Transactions on Neural
  Networks and Learning Systems}, vol.~28, no.~11, pp. 2516--2527, Nov. 2017.

\bibitem{Wang2018}
X.~Wang, H.~Wang, J.~Huang, and J.~Kurths, ``Sampled-data consensus of
  multi-agent system in the presence of packet losses,'' \emph{{IEEE} Access},
  vol.~6, pp. 54\,844--54\,853, 2018.

\bibitem{Xu2020}
L.~Xu, Y.~Mo, and L.~Xie, ``Distributed consensus over {Markovian} packet loss
  channels,'' \emph{{IEEE} Transactions on Automatic Control}, vol.~65, no.~1,
  pp. 279--286, Jan. 2020.

\bibitem{Patterson2010}
S.~Patterson and B.~Bamieh, ``Convergence rates of consensus algorithms in
  stochastic networks,'' in \emph{49th {IEEE} Conference on Decision and
  Control}.\hskip 1em plus 0.5em minus 0.4em\relax {IEEE}, Dec. 2010.

\bibitem{Wu2012}
J.~Wu and Y.~Shi, ``Average consensus in multi-agent systems with time-varying
  delays and packet losses,'' in \emph{American Control Conference}.\hskip 1em
  plus 0.5em minus 0.4em\relax {IEEE}, Jun. 2012.

\bibitem{Ghadami2012}
R.~Ghadami, ``Distributed control of multi-agent systems with switching
  topology, delay, and link failure,'' Ph.D. dissertation, Northeastern
  University, Aug. 2012.

\bibitem{Zhang2012}
Y.~Zhang and Y.-P. Tian, ``Maximum allowable loss probability for consensus of
  multi-agent systems over random weighted lossy networks,'' \emph{{IEEE}
  Transactions on Automatic Control}, vol.~57, no.~8, pp. 2127--2132, Aug.
  2012.

\bibitem{Stoorvogel2019}
A.~A. Stoorvogel, A.~Saberi, Z.~Liu, and D.~Nojavanzadeh, ``{H2} and
  {H}$\infty$ almost output synchronization of heterogeneous continuous-time
  multi-agent systems with passive agents and partial-state coupling via static
  protocol,'' \emph{International Journal of Robust and Nonlinear Control},
  vol.~29, no.~17, pp. 6244--6255, Aug. 2019.

\bibitem{Raza2022}
A.~Raza, M.~Iqbal, J.~Moon, and S.-I. Azuma, ``Performance measure of
  hierarchical structures for multi-agent systems,'' \emph{International
  Journal of Control, Automation and Systems}, vol.~20, no.~3, pp. 780--788,
  Mar. 2022.

\bibitem{Costa1997}
O.~L. d.~V. Costa, J.~B. R.~D. Val, and J.~C. Geromel, ``A convex programming
  approach to {H2} control of discrete-time {Markovian} jump linear systems,''
  \emph{International Journal of Control}, vol.~66, no.~4, pp. 557--580, Jan.
  1997.

\bibitem{Fioravanti2008}
A.~R. Fioravanti, A.~P.~C. Gon{\c{c}}alves, and J.~C. Geromel, ``{H2} filtering
  of discrete-time {Markov} jump linear systems through linear matrix
  inequalities,'' \emph{International Journal of Control}, vol.~81, pp.
  1221--1231, Jun. 2008.

\bibitem{Lee2015}
K.~Lee and R.~Bhattacharya, ``Stability analysis of large-scale distributed
  networked control systems with random communication delays: A switched system
  approach,'' \emph{Systems \& Control Letters}, vol.~85, pp. 77--83, Nov.
  2015.

\bibitem{Costa2005}
O.~L. d.~V. Costa, R.~P. Marques, and M.~D. Fragoso, \emph{Discrete-Time
  {Markov} Jump Linear Systems}.\hskip 1em plus 0.5em minus 0.4em\relax
  Springer London, 2005.

\bibitem{Costa1993}
O.~L. d.~V. Costa and M.~D. Fragoso, ``Stability results for discrete-time
  linear systems with {Markovian} jumping parameters,'' \emph{Journal of
  Mathematical Analysis and Applications}, vol. 179, no.~1, pp. 154--178, Oct.
  1993.

\bibitem{Fioravanti2013}
A.~R. Fioravanti, A.~P.~C. Gon{\c{c}}alves, and J.~C. Geromel, ``Optimal and
  mode-independent filters for generalised {Bernoulli} jump systems,''
  \emph{International Journal of Systems Science}, vol.~46, no.~3, pp.
  405--417, Jun. 2013.

\bibitem{Fioravanti2012}
A.~R. Fioravanti, A.~P.~C. Gon{\c{c}}alves, G.~S. Deaecto, and J.~C. Geromel,
  ``Equivalent {LMI} constraints: Applications to discrete-time {MJLS} and
  switched systems,'' in \emph{51st {IEEE} Conference on Decision and
  Control}.\hskip 1em plus 0.5em minus 0.4em\relax {IEEE}, Dec. 2012.

\bibitem{Schenato2007}
L.~Schenato, B.~Sinopoli, M.~Franceschetti, K.~Poolla, and S.~S. Sastry,
  ``Foundations of control and estimation over lossy networks,''
  \emph{Proceedings of the {IEEE}}, vol.~95, no.~1, pp. 163--187, Jan. 2007.

\bibitem{Horn2012}
R.~A. Horn and C.~R. Johnson, \emph{Matrix Analysis}.\hskip 1em plus 0.5em
  minus 0.4em\relax Cambridge, Oct. 2012.

\bibitem{Wu1995}
C.~W. Wu and L.~Chua, ``Synchronization in an array of linearly coupled
  dynamical systems,'' \emph{{IEEE} Transactions on Circuits and Systems},
  vol.~42, no.~8, pp. 430--447, Aug. 1995.

\bibitem{OlfatiSaber2007}
R.~Olfati-Saber, J.~A. Fax, and R.~M. Murray, ``Consensus and cooperation in
  networked multi-agent systems,'' \emph{Proceedings of the {IEEE}}, vol.~95,
  no.~1, pp. 215--233, Jan. 2007.

\bibitem{Loefberg2004}
J.~L{\"{o}}fberg, ``{YALMIP}: A toolbox for modeling and optimization in
  {MATLAB},'' in \emph{{IEEE} International Conference on Robotics and
  Automation}.\hskip 1em plus 0.5em minus 0.4em\relax {IEEE}, 2004.

\bibitem{SourceCode}
\BIBentryALTinterwordspacing
C.~Hespe, H.~Saadabadi, A.~Datar, H.~Werner, and Y.~Tang, ``Code for paper:
  Decomposition approach to multi-agent systems with {Bernoulli} packet loss,''
  Aug. 2022. [Online]. Available: \url{https://doi.org/10.5281/zenodo.7034465}
\BIBentrySTDinterwordspacing

\bibitem{Caverly2019}
R.~J. Caverly and J.~R. Forbes, ``{LMI} properties and applications in systems,
  stability, and control theory,'' Apr. 2021.

\end{thebibliography}
